\documentclass[12pt]{article}
\usepackage[french]{babel}
\usepackage{amssymb}
\usepackage{eufrak}
\usepackage{amsmath}
\usepackage{t1enc}
\usepackage{mathrsfs}
\usepackage[latin1]{inputenc}

\include{xy}
\xyoption{all}

\begin{document}
\def\QED{\hbox{\hskip 1pt \vrule width4pt height 6pt depth 1.5pt \hskip 1pt}}
\newcounter{defcount}[section]
\setlength{\parskip}{1ex}
\newtheorem{theorem}{Th{\'e}or{\`e}me}[section]
\newtheorem{lemme}[theorem]{Lemme}
\newtheorem{il}[theorem]{Illustration}
\newtheorem{coro}[theorem]{Corollaire}
\newtheorem{prop}[theorem]{Proposition}
\newtheorem{definition}[theorem]{D{\'e}finition}
\newtheorem{remark}[theorem]{Remarque}
\newtheorem{remarks}[theorem]{Remarks}
\newtheorem{ex}[theorem]{Exemples}
\newtheorem{hyp}[theorem]{Hypoth{\`e}se}
\newtheorem{Not}[theorem]{Notations}

\newcommand{\Ct}{C_{\tau}}
\newcommand{ \Kt}{K_{\tau}}
\newcommand{\st}{\sigma_{\tau}}
\newcommand{\Gt}{Aut(C_{\tau}  / \mathbb{C})}
\def\GL{{\rm GL}}
\def\SL{{\rm SL}}
\def\Sp{{\rm Sp}}
\def\sl{{\rm sl}}
\def\sp{{\rm sp}}
\def\gl{{\rm gl}}
\def\PSL{{\rm PSL}}
\def\SO{{\rm SO}}
\def\Gal{{\rm Gal}}
\def\Sym{{\rm Sym}}
\def\tildea{\tilde a}
\def\tildeb{\tilde b}
\def\tildeg{\tilde g}
\def\tildee{\tilde e}
\def\tildeA{\tilde A}
\def\calG{{\cal G}}
\def\calH{{\cal H}}
\def\calT{{\cal T}}
\def\calC{{\cal C}}
\def\calP{{\cal P}}
\def\calS{{\cal S}}
\def\calM{{\cal M}}
\def\curve{{\rm {\bf C}}}
\def\P1{{\rm {\bf P}}^1}
\def\Ad{{\rm Ad}}
\def\ad{{\rm ad}}
\def\Aut{{\rm Aut}}
\def\Int{{\rm Int}}
\def\C{{\mathbb C}}
\def\QX{{\mathbb Q}}
\def\NX{{\mathbb N}}
\def\ZX{{\mathbb Z}}
\def\DX{{\mathbb D}}
\def\AX{{\mathbb A}}
\def\semi{\hbox{${\vrule height 5.2pt depth .3pt}\kern -1.7pt\times $ }}

\newenvironment{prf}[1]{\trivlist
\item[\hskip \labelsep{\bf
#1.\hspace*{.3em}}]}{~\hspace{\fill}~$\square$\endtrivlist}
\newenvironment{proof}{\begin{prf}{Preuve}}{\end{prf}}
 \def\square{\QED}
 \newenvironment{proofofthm}{\begin{prf}{Proof of Theorem~\ref{ext}}}{\end{prf}}
 \def\square{\QED}
\newenvironment{sketchproof}{\begin{prf}{Sketch of Proof}}{\end{prf}}

\begin{titlepage}
\begin{center}
\vspace{2.8cm} {\LARGE\bf Hypertranscendance  }
\end{center}

\begin{center}
{\LARGE\bf   et Groupes de Galois
aux diff{\'e}rences}
\end{center}\vspace{.9cm}

\begin{center}
  {\large{{\bf  Charlotte
Hardouin}}}\\\vspace{.5cm} {\large{{\bf Universit{\'e} Paris VI}}}
\end{center}

\begin{center} le 20  septembre   2006 \\
\end{center}
\vspace{.9cm}
\end{titlepage}

{\bf R{\'e}sum{\'e}}

\medskip
 On donne  dans cet article divers crit{\`e}res
 d'ind{\'e}pendance alg{\'e}brique pour les d{\'e}riv{\'e}es successives de solutions d'{\'e}quations aux diff{\'e}rences de rang 1.  L'id{\'e}e principale consiste {\`a} construire au moyen de l'op{\'e}rateur de d{\'e}rivation, qui commute avec l'op{\'e}rateur aux diff{\'e}rences, des extensions it{\'e}r{\'e}es du module aux diff{\'e}rences initial. Le probl{\`e}me se ram{\`e}ne alors  au calcul du groupe de Galois aux diff{\'e}rences de telles extensions, calcul qui lui-m{\^e}me se r{\'e}duit {\`a} une simple question d'alg{\`e}bre lin{\'e}aire. Les cat{\'e}gories tannakiennes mises en jeu sont neutres, mais sur des corps parfois non alg{\'e}briquement clos,
  ce qui conduit {\`a}  {\'e}tudier le comportement des groupes de Galois par extension des corps de bases.

\bigskip
{\bf Abstract}

\medskip
 This paper deals with criteria of algebraic
independence for the derivatives of solutions of rank one difference
equations. The key idea consists in deriving from the commutativity
of the differentiation and difference operators a sequence of
iterated extensions of the original difference module, thereby
setting the problem in the framework of difference Galois theory and
finally reducing it to an exercise in linear algebra. The involved
tannakian categories are neutral over non necessarily algebraically
closed fields, and this leads us to study the behaviour of Galois
groups under base field extensions.

\bigskip

\vspace{2cm}

- Key words: linear difference equations, difference Galois groups,
  tannakian category, hypertranscendence, Galois cohomology

- Classification AMS: 12 H 10, 12 H 05,  34 M 15, 39 A 13, 11 R 34

-  Adresse electronique de l'auteur:  hardouin@math.jussieu.fr,

- Adresse postale : \begin{enumerate}
\item[.]Institut de
math{\'e}matiques de Jussieu, 175 rue du Chevaleret, 75013 Paris
\item[.] {\`a} compter du 01/10/2006 :  Interdisziplin{\"a}res Zentrum f{\"u}r Wissenschaftliches Rechnen
(IWR)  der Universit{\"a}t Heidelberg, Im Neuenheimer Feld 368,
69120 Heidelberg Germany
\end{enumerate}
\bigskip

\vfil \eject

\tableofcontents

\section{Introduction}

Dans cet article, on s'int{\'e}resse aux relations alg{\'e}bro-diff{\'e}rentielles
satisfaites par les  solutions d'{\'e}quations fonctionnelles, plus particuli{\`e}rement aux diff{\'e}rences. On dira qu'une
fonction $f$ (ind{\'e}finiment d{\'e}rivable) est \textit{hypertranscendante} sur un corps  $K$ s'il n'existe pas de
relations alg{\'e}briques   {\`a} coefficients dans $K$ liant $f$ et ses
d{\'e}riv{\'e}es, autrement dit si $f$ ne v{\'e}rifie pas d'{\'e}quation diff{\'e}rentielle alg{\'e}brique sur $K$.\\
L'exemple le plus classique de fonction hypertranscendante est celui
de la fonction $\Gamma$ qui v{\'e}rifie l'{\'e}quation fonctionnelle
$\Gamma(z+1)= z \Gamma(z)$ et qui est hypertranscendante sur
$\mathbb{C}(z)$.\\
Des r{\'e}sultats de m{\^e}me nature ont {\'e}t{\'e} {\'e}tablis par exemple  par S. Bank (voir \cite{Bank}) dans le cadre
d'{\'e}quations aux $\tau$-diff{\'e}rences de  rang $1$ et par K. Ishizaki (\cite{Ishi}) dans
le cadre d'{\'e}quations  aux $q$-diff{\'e}rences non homog{\`e}nes de rang
$1$. \\
Ces r{\'e}sultats se basent essentiellement sur  des m{\'e}thodes analytiques
et sur  des estimations de la taille des coefficients du d{\'e}veloppement de
telles fonctions, telles qu'{\'e}tudi{\'e}es dans \cite{Ram}.\\
Dans cet article, on cherche à {\'e}tendre ces {\'e}nonc{\'e}s {\`a}
l'{\'e}tude simultan{\'e}e de plusieurs solutions (autrement dit,
{\`a} des questions d ' \textit{"hyperind{\'e}pendance
alg{\'e}brique"}), et ce, par des m{\'e}thodes purement
alg{\'e}briques de th{\'e}orie de Galois aux diff{\'e}rences. Voici,
dans le cadre des $q$-diff{\'e}rences, le type de
th{\'e}or{\`e}me auquel on aboutira.\\

Soit $q \in \mathbb{C}^*, |q| \neq 1$. On d{\'e}signe par
$\mathcal{M}er(\mathbb{C}^*)$ le corps des fonctions m{\'e}romorphes sur
$\mathbb{C}^*$, par $\sigma_q$ l'automorphisme de
$\mathcal{M}er(\mathbb{C}^*)$, qui {\`a} une fonction $f(z)$ associe la
fonction $f(qz)$ et par $C_E$ le sous-corps de
$\mathcal{M}er(\mathbb{C}^*)$ des fonctions invariantes sous l'action
de $\sigma_q$, qui est isomorphe au  corps $\mathbb{C}(E)$ des fonctions rationnelles
sur la courbe elliptique   $E=\mathbb{C}^* / q^{\mathbb{Z}}$.\\
On appelle diviseur elliptique d'une fonction $f \in \mathbb{C}(z)$
l'image dans le groupe des diviseurs de la courbe elliptique
$E = \mathbb{C}^* / q^{\mathbb{Z}}$ de la partie premi{\`e}re {\`a} $0$ du
diviseur de $f$.

\begin{theorem}\label{theorem:gen}
Soient $a_1,...,a_n$ des {\'e}l{\'e}ments non nuls de $\mathbb{C}(z)$ et $q$
un nombre complexe non nul de module diff{\'e}rent de $1$. Pour tout
$i=1,..., n$, soit $f_i \neq 0$  une  solution m{\'e}romorphe sur
$\mathbb{C}^*$ de l'{\'e}quation aux $q$-diff{\'e}rences
$\sigma_q(f_i)=a_i f_i$. On suppose que les diviseurs elliptiques
des $a_i$ sont lin{\'e}airement ind{\'e}pendants sur $\mathbb{Z}$.
Alors les fonctions $f_1,...,f_n$ ainsi que leur d{\'e}riv{\'e}es
successives sont alg{\'e}briquement ind{\'e}pendantes sur $C_E(z)$.
\end{theorem}
On obtient un {\'e}nonc{\'e} similaire pour les $\tau$-diff{\'e}rences, qui
entra{\^\i}ne, par exemple, que les seules relations de d{\'e}pendance
alg{\'e}briques liant les fonctions polygamma $\psi^{(k)}(z
+\alpha),\psi^{(k)}(nz), k, n
\in \mathbb{N}, \alpha \in \mathbb{C}$ se d{\'e}duisent de la relation de
distribution satisfaite par la fonction gamma.\\

L'outil essentiel de  la d{\'e}monstration du th{\'e}or{\`e}me $1.1$  est la th{\'e}orie de
Galois aux diff{\'e}rences, qui permet de faire le lien entre groupe de
Galois des {\'e}quations et degr{\'e} de transcendance des extensions de corps
mises en jeu (tout comme en th{\'e}orie  des {\'e}quations diff{\'e}rentielles). Le principe de la preuve est le suivant.

Soit $(K,\sigma,\partial)$ un corps aux diff{\'e}rences  et diff{\'e}rentiel
 de caract{\'e}ristique nulle, c'est {\`a} dire un corps $K$ muni d'un
 automorphisme de corps $\sigma$ et d'un op{\'e}rateur de d{\'e}rivation
 $\partial$ tels que  $\sigma \circ \partial =\partial \circ
\sigma$. On note $C_{\sigma}$ (resp. $C_{\partial}$) le sous-corps de $K$ form{\'e}  des
 invariants sous $\sigma$ (resp.  des
 constantes diff{\'e}rentielles  de $K$) et  $\mathcal{D}_{\sigma}
 =K[\sigma,\sigma^{-1}]$ l'anneau des polyn{\^o}mes en $\sigma$ et
 $\sigma^{-1}$ {\`a} coefficients dans $K$. Par exemple dans le cas supra, on  pourra prendre
 $(K,\sigma,\partial)=(C_E(z),\sigma_q, zd/dz)$
 (d'o{\`u} $C_{\sigma}=C_E$ et $C_{\partial}=\mathbb{C}$), ou encore $(K,
 \sigma, \partial)=(\mathbb{C}(z),\sigma_q, zd/dz)$(d'o{\`u} $C_{\sigma}=C_{\partial}=\mathbb{C}$).

On consid{\`e}re le syst{\`e}me  aux diff{\'e}rences :
\begin{equation} \label{eqn:der1q}
\sigma Y=AY, \ \mbox{o{\`u}} \ A   \in Gl_n(K) \end{equation}

Soient $L$ une extension aux diff{\'e}rences diff{\'e}rentielle de $K$ et $Y \in
 {L}^n$ une solution de (\ref{eqn:der1q}).
 Notons $\mathcal{A}$ le $\mathcal{D}_{\sigma}$-module associ{\'e} {\`a} ce
 syst{\`e}me et consid{\'e}rons le vecteur  $$  \left(\begin{array}{c} \partial Y \\
                                       Y \end{array} \right)\in L^{2n}.$$
                                   Puisque $\sigma$ et $\partial$
 commutent, il   v{\'e}rifie l'{\'e}quation
                                     aux diff{\'e}rences :
\begin{equation} \label{eqn:der2q} \sigma \left(\begin{array}{c} \partial Y \\
                                       Y \end{array} \right)
                                  = \left(\begin{array}{cc} A & \partial A \\
                                    0   & A \end{array} \right)
                                     \left(\begin{array}{c} \partial Y \\
                                       Y \end{array} \right) \
                                       \end{equation}
qui correspond {\`a} une extension $\mathcal{M}(1)$ de $\mathcal{A}$ par
$ \mathcal{A}$
                         dans la cat{\'e}gorie des
                                       $\mathcal{D}_{\sigma}$-modules.\\

En it{\'e}rant $m-1$ fois ce processus, et en notant que pour tout entier $j \geq 1$,
\begin{equation}\label{eqn:iterder}  \sigma(\partial^j Y)= \partial^j(\sigma Y) =
\partial^j(AY)= \sum_{i=0}^j C^i_j \partial^i A
\partial^{j-i}Y, \end{equation}
 on obtient le   syst{\`e}me
 aux   diff{\'e}rences\\
 \begin{equation} \label{eqn:dern} \sigma   \left(
 \begin{array}{c}
\partial^mY \\

.\\
.\\
.\\
.\\
.\\
.\\
.\\
\partial Y\\
                 Y \end{array}\right)=    \left(\begin{array}{ccccccc}\\
        A & \cdots & \cdots & C_m^{k} \partial^{k}A & \cdots  &C_m^1 \partial^{m-1} A & \partial^{m} A   \\
        0 & A  & \cdots & \vdots    & \cdots  & \cdots   & \vdots\\
        0 &\cdots &  \ddots  &    C^{k-r}_{m-r}\partial^{k-r}A  &... &\cdots   & \partial^{m-r}A \\
        0   & ... &... & \cdots & \cdots & \cdots & \vdots\\
        0   & ... &... & A & \cdots & \cdots & \vdots\\
   0&...& ...& ... & \ddots&  \cdots & \partial^2 A \\
     0&...& ...& ... &...&  A & \partial A \\
         0& ...&...&...&...& 0 & A \end{array} \right)  \left(
 \begin{array}{c}
\partial^mY \\

.\\
.\\
.\\
.\\
.\\
.\\
.\\
\partial Y\\
                 Y \end{array}\right); \end{equation} c'est la
                 repr{\'e}sentation matricielle  d'un $
                 \mathcal{D}_{\sigma}$-module $\mathcal{M}(m)$,  extension it{\'e}r{\'e}e
                 $m$-fois de l'objet $\mathcal{A}$
                 dans cette cat{\'e}gorie.\\

Le degr{\'e} de transcendance du corps $K(Y,\partial
Y,...,\partial^m Y)$ sur $K$ s'interpr{\`e}te alors comme la
dimension de l'orbite de cette solution sous l'action du groupe de
Galois aux diff{\'e}rences  $G(\mathcal{M}(m))$ de $\mathcal{M}(m)$.
De plus, lorsque $\mathcal{A}$ est compl{\'e}tement r{\'e}ductible
la seule action du radical unipotent de  $G(\mathcal{M}(m))$ suffit
{\`a} contr{\^o}ler le degr{\'e} de transcendance de la \textit{partie
diff{\'e}rentielle}  $K(Y,
\partial
Y,...,\partial^m Y) / K(Y)$ de cette extension. \\
C'est ainsi  au calcul du  radical unipotent de $G(\mathcal{M}(m))$
qu'est consacr{\'e} l'essentiel de l'article, pour $\mathcal{A}$
somme directe de $\mathcal{D}_{\sigma}$-modules de rang $1$. Sous
les hypoth{\`e}ses du th{\'e}or{\`e}me \ref{theorem:gen}, on verra
qu'il est aussi grand que possible. Le th{\'e}or{\`e}me
\ref{theorem:gen} et son analogue pour les
$\tau$-diff{\'e}rences en d{\'e}coulent imm{\'e}diatement.\\

Le plan de l' article est le suivant.\\

 Dans la partie $2$, on {\'e}tend  les
 r{\'e}sultats de \cite{BS} et  \cite{Ber}  sur  le calcul
 du radical unipotent du  groupe de Galois du produit de deux
 op{\'e}rateurs diff{\'e}rentiels  compl{\'e}tement r{\'e}ductibles  au cadre d'une
 cat{\'e}gorie tannakienne $\bold{T}$ neutre sur un corps $C$ de
 caract{\'e}ristique $0$ mais non n{\'e}cessairement alg{\'e}briquement clos. Plus
 pr{\'e}cis{\'e}ment, on d{\'e}termine au th{\'e}or{\`e}me \ref{theorem:Ber2}
 le radical unipotent du groupe de Galois d'une extension  de l'objet $\bold{1}$, {\'e}l{\'e}ment unit{\'e} de $\bold{T}$,
 par un objet $\mathcal{Y}$ compl{\'e}tement r{\'e}ductible. On obtient ainsi
 l'{\'e}quivalence entre la complexit{\'e} de l'extension et la taille de son
 groupe de Galois. On rappelle {\'e}galement (proposition \ref{prop:degtr}) comment le
 groupe de Galois contr{\^o}le le degr{\'e} de transcendance des solutions.\\

On peut extraire de  la litt{\'e}rature  plusieurs d{\'e}finitions de groupes de
Galois aux diff{\'e}rences, correspondant {\`a} divers choix de foncteurs fibres. Par exemple, dans le cas des $q$-diff{\'e}rences,
celui des symboles  de M. van der  Put et M.F. Singer (\cite{VPS2}) o{\`u} la cat{\'e}gorie tannakienne $\bold{T}$
est la cat{\'e}gorie des $\sigma_q$-modules de type
fini sur $K=\mathbb{C}(z)$ (ou plus g{\'e}n{\'e}ralement sur $K=C(z)$ avec
$C$ alg{\'e}briquement clos) et  o{\`u} le groupe de
Galois correspond au choix d'un foncteur fibre $\omega$ sur $\mathbb{C}$
(voir aussi \cite{And},  paragraphes $ 3.2$ et $3.4$). Les foncteurs fibres, de nature
plus g{\'e}om{\'e}triques, de J. Sauloy  (\cite{Sau}) fournissent encore des
groupes de Galois sur $\mathbb{C}$.   Les travaux de C. Praagman (\cite{Prag})
permettent en revanche  de consid{\'e}rer la cat{\'e}gorie
 $\bold{T}_E$ des $\sigma_q$-modules de type
fini sur $K_E=C_E(z)$ et un foncteur fibre   $\omega_E$  (donc un groupe de Galois) sur le corps non alg{\'e}briquement clos  $C_E$.
La nature m{\^e}me de nos r{\'e}sultats ({\'e}tude de ``vraies" fonctions, sur lesquelles agit une d{\'e}rivation) conduit {\`a}  privil{\'e}gier cette derni{\`e}re approche \footnote{ On
  trouvera au paragraphe 3.3 une d{\'e}monstration fond{\'e}e sur la premi{\`e}re
  approche et sur un cas particulier des th{\'e}or{\`e}mes de comparaison
  {\'e}tablis dans \cite{CHS}   entre groupes de Galois au sens de
  $\bold{T}$ et de $\bold{T}_E$}, mais on montre dans la partie  $3$ que
des  extensions d'objets de $\bold{T}$ ind{\'e}pendantes au sens de
$\bold{T}$ le demeurent, apr{\`e}s  extension des scalaires {\`a} $K_E$, au
sens de $\bold{T}_E$ (th{\'e}or{\`e}me \ref{theorem:injext}). On {\'e}tablit
d'ailleurs ce fait dans un cadre plus g{\'e}n{\'e}ral, inspir{\'e} des m{\'e}thodes de
cohomologie galoisienne.

Dans la partie  $4$, on v{\'e}rifie que les extensions de la cat{\'e}gorie
$\bold{T}$ mises en jeu au th{\'e}or{\`e}me \ref{theorem:gen} sont bien
ind{\'e}pendantes. Dans le cadre particulier des {\'e}quations aux
$q$-diff{\'e}rences de rang $1$ ou  plus g{\'e}n{\'e}ralement des syst{\`e}mes
diagonaux, il s'agit l{\`a} d'une simple question d'alg{\`e}bre lin{\'e}aire (cf
remarque $4.14$). On aboutit ainsi aux r{\'e}sultats d'hypertranscendance
du th{\'e}or{\`e}me  \ref{theorem:gen}.\\

Enfin, on donne dans la partie  $5$ l'analogue des th{\'e}or{\`e}mes
de la partie  $4$ dans le cadre des $\tau$-diff{\'e}rences.

\textbf{Remarque}\\
La consid{\'e}ration simultan{\'e}e de deux op{\'e}rateurs, ici
$\sigma_q$ et $\partial$, agissant sur le m{\^e}me corps de
fonctions n'est  pas  nouvelle, voir par exemple les travaux de J.P.
B{\'e}zivin (\cite{Bez2}) pour deux op{\'e}rateurs aux
$q$-diff{\'e}rences, $\sigma_{q_1}$ et $\sigma_{q_2}$ avec $q_1 \neq
q_2$. Dans ce type de travail, les deux op{\'e}rateurs jouent un r{\^o}le sym{\'e}trique et on ne traite que d'{\'e}quations lin{\'e}aires en chacun des op{\'e}rateurs.\\
L'id{\'e}e cl{\'e} du pr{\'e}sent article est la consid{\'e}ration des extensions
successives $\mathcal{M}(n)$ d{\'e}crites plus haut. Les r{\^o}les des op{\'e}rateurs
$\sigma_q$ et $\partial$ ne sont plus sym{\'e}triques, permettant ainsi
d'atteindre des {\'e}quations diff{\'e}rentielles non lin{\'e}aires. Signalons
que ce processus s'adapte {\`a} d'autres situations : voir ainsi le
recent travail de A. Ovchinnikov \cite{Ov}, o{\`u} $\sigma_q$ est
remplac{\'e} par un operateur de
derivation $\partial_x$ commutant avec $\partial_z$.\\

\textbf{Remerciements}

Je tiens a remercier sincérement D.Bertrand de m'avoir menée de
façon si attentive sur le chemin de ces questions
d'hypertranscendence et de m'avoir aussi bien assistée à chaque pas
du chemin. Je voudrais aussi remercier M. Singer et Z. Chatzidatkis
pour nos discussions communes et pour notre travail sur la
définition des groupes de Galois aux différences.

\section{Calcul de groupes de Galois  }

Dans cette partie, $C$ d{\'e}signe un corps de caract{\'e}ristique
nulle. {\it On ne suppose  pas que $C$ soit alg{\'e}briquement
clos.}\\
 Soit  $\bold{T}$ une cat{\'e}gorie tannakienne neutre sur
le corps
 $C$ (au sens de  Deligne \cite{Dlct}, dont on reprend dans tout l'article les d{\'e}finitions).  On note $\bold{1}$ l'objet
unit{\'e} (en particulier $End(\bold{1})=C$) et $\omega$ un foncteur
fibre de $\bold{T}$ {\`a} valeurs dans la cat{\'e}gorie $Vect_C$ des
 $C$-espaces vectoriels de dimension finie.\\

 Le principal r{\'e}sultat de cette partie est la  version tannakienne
 suivante   des {\'e}nonc{\'e}s de \cite{BS} et \cite{Ber}, o{\`u} l'on  appelle \textit{groupe de
Galois} d'un objet $\mathcal{M}$ de $\bold{T}$, not{\'e}
$G_{\mathcal{M}}$, le groupe alg{\'e}brique sur $C$ repr{\'e}sentant le foncteur
$Aut^{\otimes}(\omega|<\mathcal{M}>)$ (o{\`u} $<\mathcal{M}>$
d{\'e}signe la cat{\'e}gorie tannakienne engendr{\'e}e par
$\mathcal{M}$ dans $\bold{T}$).
\medskip

\begin{theorem}\label{theorem:Ber2}
Soient
$\mathcal{Y}$ un objet compl{\'e}tement r{\'e}ductible  de
$\bold{T}$  et $\mathcal{U}$ une extension dans $\bold{T}$ de
$\bold{1}$ par  $\mathcal{Y}$. Alors $G_{\mathcal{U}}$ est
{\'e}gal au produit semi-direct du groupe r{\'e}ductif  $G_{\mathcal{Y}}$ par le  groupe vectoriel $\omega(\mathcal{V})$,
o{\`u} $\mathcal{V}$ d{\'e}signe le plus petit sous-objet de
$\mathcal{Y}$ tel que le quotient par
$\mathcal{V}$ de l'{\'e}l{\'e}ment  $\mathcal{U}$ de
$Ext^1(\bold{1},\mathcal{Y})$ soit une
extension scind{\'e}e.\\
\end{theorem}

L'existence du plus petit objet $\mathcal{V}$ de l'{\'e}nonc{\'e} pr{\'e}c{\'e}dent
est justifi{\'e} au lemme \ref{lemme:pp} et la preuve du th{\'e}or{\`e}me est
donn{\'e}e au paragraphe $2.2$.

On d{\'e}duit du th{\'e}or{\`e}me  \ref{theorem:Ber2} le corollaire suivant, qui suffira pour
les applications d{\'e}crites aux parties $4$ et $5$ de cet article.\\

\begin{coro}\label{coro:inde}
 Soient
  $\mathcal{Y}$ un  objet compl{\'e}tement r{\'e}ductibles  de $\bold{T}$, $\Delta$ l'anneau
  $End(\mathcal{Y})$,
  $ \mathcal{E}_1,...,\mathcal{E}_n$, des   extensions  de $\bold{1}$
  par $ \mathcal{Y} $   telles que
  $\mathcal{E}_1,...,\mathcal{E}_n$ soient
  $\Delta$-lin{\'e}airement ind{\'e}pendantes dans $Ext^1_{\bold{T}}(\bold{1},\mathcal{Y})$. Alors le
radical unipotent de $G_{\mathcal{E}_1 \oplus... \oplus
\mathcal{E}_n}$ est isomorphe {\`a} $ \omega(\mathcal{Y})^n$.\\
\end{coro}

On conclut cette seconde partie   par le lien entre le degr{\'e} de transcendance des solutions
et la dimension du groupe de Galois mis en jeu.\\

\subsection{Groupe de Galois d'objets compl{\'e}tement r{\'e}ductibles}

Soient donc  $C$ un corps  commutatif de caract{\'e}ristique nulle,
 $\bold{T}$ une cat{\'e}gorie tannakienne
neutre  sur $C$, $\bold{1}$ l'objet neutre. On fixe d{\'e}sormais $\omega$ un foncteur
fibre de $\bold{T}$ dans la cat{\'e}gorie $Vect_C$ des espaces
vectoriels de dimension
finie et on s'autorise {\`a} omettre l'indice $\omega$ dans les
 d{\'e}finitions qui suivent.\\
 On d{\'e}signe par $\mathcal{G}=  {\mathcal{G}}^{\omega}  $ le groupe proalg{\'e}brique sur
$C$ repr{\'e}sentant le foncteur $Aut^{\otimes}(\omega)$, par
$Rep_{\mathcal{G}}$ la cat{\'e}gorie des repr{\'e}sentations de
$\mathcal{G}$ port{\'e}es par des $C$-espaces vectoriels de
dimension finie sur $C$, et par  $1_C$ l'objet neutre de
$Rep_{\mathcal{G}}$. On peut d{\'e}crire les quotients de
$\mathcal{G}$ de type $G_{\mathcal{M}} =
Aut^{\otimes}(\omega|<\mathcal{M}>)$ de la fa\c {c}on suivante.\\

\begin{definition}
Pour tout objet $\mathcal{M}$  de $\bold{T}$, on note
$Stab(\mathcal{M})$ le sous-sch{\'e}ma en  groupe de $\mathcal{G}$
v{\'e}rifiant : pour tout $ \sigma \in Stab(\cal M),$ et  tout  $  x
\ \in \omega(\mathcal{M}),\  \sigma(x)= x$. Alors,
$Stab(\mathcal{M})$ est un sous-sch{\'e}ma en groupes  distingu{\'e}
de  $\mathcal{G}$, et $\rm{Gal}(\mathcal{M})=
{\mathcal{G}}/Stab({\cal M}) \simeq$ $ G_{\mathcal{M}} $ est
naturellement muni d'une structure de groupe alg{\'e}brique sur $C$.
 On dit que $\rm{Gal}(\mathcal{M})$ est le \textit{groupe
 de Galois} de $\mathcal{M}$ relatif {\`a} $\omega$.
\end{definition}

\begin{definition} Soit $\mathcal{X}$ un objet de  $\bold{T}$. On dit que $\mathcal{X}$
 est \textit{irr{\'e}ductible} s'il n'admet pas de sous-objet propre.
On dit que $\mathcal{X}$ est \textit{isotypique} s'il est somme
directe d'objets irr{\'e}ductibles isomorphes.  On dit que $\mathcal{X}$
 est \textit{compl{\'e}tement r{\'e}ductible} s'il est somme directe
 d'objets irr{\'e}ductibles ou de fa{\c c}on {\'e}quivalente s'il est somme directe
 de ses composantes isotypiques maximales.
\end{definition}

\begin{prop}\label{prop:reduc}
Soit $\mathcal{X} \in \bold{T}$ et   $Rep_\mathcal{G}$  comme ci-dessus.\\
 Les propri{\'e}t{\'e}s suivantes sont
{\'e}quivalentes.
\begin{enumerate}
\item $\mathcal{X}$ est \textit{compl{\'e}tement r{\'e}ductible} dans $  \bold{T}$.
\item $\omega(\mathcal{X})$ est \textit{compl{\'e}tement r{\'e}ductible} dans
  $Rep_\mathcal{G}$.
\item $\rm{Gal}(\mathcal{X})$ est un groupe r{\'e}ductif, c'est-{\`a}-dire son
  radical unipotent est trivial.
\item Tout $\rm{Gal}(\mathcal{X})$-module est compl{\'e}tement r{\'e}ductible.
\end{enumerate}
\end{prop}

\textbf{D{\'e}monstration}\\
Le foncteur fibre $\omega$ fournit une {\'e}quivalence de
cat{\'e}gorie entre $\bold{T}$ et  $Rep_\mathcal{G}$ (\cite{Dlct}).
L'{\'e}quivalence entre $1)$ et $2)$ en d{\'e}coule. Les
{\'e}nonc{\'e}s $2)$ et $4)$ sont {\'e}quivalents car  $\mathcal{G}$
agit sur $\omega(\mathcal{X})$ {\`a} travers $Gal(\mathcal{X})$, et
si  un groupe admet une repr{\'e}sentation fid{\`e}le (ici $\omega
(\mathcal{X})$) compl{\'e}tement r{\'e}ductible, toutes ses
repr{\'e}sentations sont alors
compl{\'e}tement r{\'e}ductibles.\\
D'apr{\`e}s \cite{DM}, proposition 2.23  et
remarque 2.28, p.141, un groupe alg{\'e}brique $G$
d{\'e}fini sur un corps quelconque  de caract{\'e}ristique nulle est r{\'e}ductif si et
seulement si la cat{\'e}gorie $Rep_G$ est semi-simple. Les deux derni{\`e}res
propri{\'e}t{\'e}s sont donc {\'e}quivalentes.


Pour tout objet  $\mathcal{Y}$  (resp. $Y$ )
d'objets de $\bold{T}$ (resp. de $Rep_{\mathcal{G}}$), notons
$Ext^1_{\bold{T}}(\bold{1},\mathcal{Y})$ (resp.
$Ext^1_{Rep_{\mathcal{G}} }(1_C,Y)$) le groupe des classes
d'isomorphismes d'extensions de $\bold{1}$ (resp. $1_C$) par
$\mathcal{Y}$ (resp. $Y$). Ces groupes sont naturellement munis d'une
action de $C$, qui en fait des $C$-espaces vectoriels. Plus
g{\'e}n{\'e}ralement, ces groupes sont naturellement munis d'une structure de
module sur la $C$-alg{\`e}bre  $ \Delta =
End_{\bold{T}}(\mathcal{Y}) \simeq  End_{Rep_{\mathcal{G}}}(Y)$).\\

\begin{prop} \label{prop:cocycle}

  Pour tout couple $(\mathcal{X}, \mathcal{Y})$ d'objets  de $\bold{T}$,  les $C$-espaces vectoriels
  $Ext^1_{\bold{T}}(\mathcal{X},\mathcal{Y})$,
$Ext^1_{Rep_\mathcal{G}}(\omega(\mathcal{X}),\omega(\mathcal{Y}))$
sont canoniquement isomorphes.
\end{prop}

\textbf{D{\'e}monstration}\\
L'isomorphisme  d{\'e}coule de l'{\'e}quivalence de cat{\'e}gories entre $\bold{T}$
et
$Rep_{\mathcal{G}}$.\\

\subsection{Groupe de Galois d'une extension simple}

On reprend les hypoth{\`e}ses et les notations  du th{\'e}or{\`e}me \ref{theorem:Ber2}.
On note de plus  $G_{\mathcal{Y}}=Gal(\mathcal{Y})$ et
$G_{\mathcal{U}}=Gal(\mathcal{U})$, et on pose
$$R_u= Stab(\mathcal{Y}) / Stab(\mathcal{U}) .$$
 Ainsi $R_u$ est distingu{\'e} dans
  $G_{\mathcal{U}}$ et   $G_\mathcal{Y}$ est isomorphe  {\`a}
$G_\mathcal{U} / R_u   $. On verra au cours de la d{\'e}monstration que
 $R_u$ est un sous-groupe vectoriel de $G_{\mathcal{U}}$. Puisque
 $\mathcal{Y}$ est, par hypoth{\`e}se, compl{\'e}tement r{\'e}ductible,
 le groupe $G_{\mathcal{U}}/R_u \simeq G_{\mathcal{Y}}$ est r{\'e}ductif
(cf. Proposition \ref{prop:reduc}), et $R_u$ est donc le radical unipotent de $G_{\mathcal{U}}$.\\

Avant de passer {\`a} la preuve du th{\'e}or{\`e}me \ref{theorem:Ber2}, on d{\'e}montre
le lemme pr{\'e}liminaire suivant.

\begin{lemme}[et d{\'e}finition]\label{lemme:pp}
Soient $\mathcal{Y}$ un  objet compl{\'e}tement r{\'e}ductible de $\bold{T}$, et
$\mathcal{U}$ un objet de $\bold{T}$ extension de $\bold{1}$ par
$\mathcal{Y}$. L'ensemble $\bold{V}= \bold{V}(\mathcal{Y},\mathcal{U})$ des  sous-objets $\mathcal{V}$ de $\mathcal{Y}$ tel que le
quotient  de l'extension $\mathcal{U}$ par
$\mathcal{V}$  est une extension  scind{\'e}e  admet un plus petit
{\'e}l{\'e}ment pour l'inclusion. On l'appellera \textbf{l'objet  minimal}  de $\bold{V}(\mathcal{Y},\mathcal{U})$.

\end{lemme}

\textbf{D{\'e}monstration}\\
Il s'agit de voir que si $\mathcal{V}_1$ et $\mathcal{V}_2$ sont des
{\'e}l{\'e}ments de $\bold{V}$, il en est de m{\^e}me de  leur
intersection $\mathcal{W}$. Puisque $\mathcal{Y}$ est
compl{\`e}tement r{\'e}ductible,   il existe trois sous-objets
$\mathcal{V'}$, $\mathcal{W}_1'$,  $\mathcal{W}_2'$  de
$\mathcal{Y}$ tels que :
\begin{enumerate}
\item $\mathcal{V}_1=\mathcal{W} \oplus \mathcal{W}_1'$, $ \mathcal{V}_2=\mathcal{W} \oplus \mathcal{W}_2'$.
\item $\mathcal{Y} = \mathcal{V}_1 \oplus \mathcal{W}_2'\oplus \mathcal{V}' = \mathcal{V}_2 \oplus \mathcal{W}_1' \oplus \mathcal{V}'=\mathcal{W} \oplus \mathcal{W}_2' \oplus \mathcal{W}_1' \oplus \mathcal{V}$
\end{enumerate}
On a: $$Ext^1(\bold{1}, \mathcal{Y})\simeq
Ext^1(\bold{1}, \mathcal{V}_1) \times Ext^1(\bold{1},\mathcal{W}_2'\oplus \mathcal{V'}) \ \mbox{ et}
\  Ext^1(\bold{1}, \mathcal{Y})\simeq
Ext^1(\bold{1}, \mathcal{V}_2) \times Ext^1(\bold{1},\mathcal{W}_1'\oplus \mathcal{V'}) .$$ Comme $\mathcal{V}_1$ et $\mathcal{V}_2$ sont dans $\bold{V}$,
$\mathcal{U}$ se projette de fa{\c c}on triviale sur
$Ext^1(\bold{1}, \mathcal{W}_2'\oplus \mathcal{V'})$ et sur  $ Ext^1(\bold{1},\mathcal{W}_1'\oplus \mathcal{V'} )$. On
en d{\'e}duit que $ \mathcal{U}$ se projette de fa{\c c}on
triviale sur $Ext^1(\bold{1},\mathcal{W}_2'\oplus \mathcal{W}_1'\oplus \mathcal{V'})$ et ainsi, que $\mathcal{W}$ est
dans $\bold{V}$.

Pour {\'e}tablir le th{\'e}or{\`e}me \ref{theorem:Ber2}, il reste donc {\`a} montrer que $R_u$ est isomorphe au groupe vectoriel $\omega(\mathcal{V})$, o{\`u}
$\mathcal{V}$ est  l'objet  minimal  de $\bold{V}(\mathcal{Y},\mathcal{U})$.\\

\textbf{D{\'e}monstration du th{\'e}or{\`e}me \ref{theorem:Ber2}}\\
Par hypoth{\`e}se, $\mathcal{U}$ s'inscrit dans une suite exacte
$$\xymatrix{
0 \ar[r] &\mathcal{Y} \ar[r]^{i} & \mathcal{U} \ar[r]^{p} & \bold{1}
\ar[r] & 0}.$$ D'apr{\`e}s la proposition \ref{prop:cocycle},
$\omega(\mathcal{U})$ se r{\'e}alise comme extension de la
repr{\'e}sentation unit{\'e} $1_C$ par $\omega(\mathcal{Y})$ dans la
cat{\'e}gorie des $\mathcal{G}$-modules de dimension finie. Soit $s$
une section $C$-lin{\'e}aire de la suite exacte de $C$-espaces
vectoriels
$$\xymatrix{
0 \ar[r] &\omega(\mathcal{Y}) \ar[r]^{\omega(i)} &
\omega(\mathcal{U}) \ar[r]^{\omega(p)} & C \ar@{.>}@/^/[l]^{s}
\ar[r] & 0}$$  \\
et posons $f=s(1) \ \in \omega(\mathcal{U})(C)$.\\
Consid{\'e}rons le morphisme  de $C$-sch{\'e}mas $ \zeta_{\omega(\cal  U)} :
G_{\mathcal{U}} \rightarrow \omega(\mathcal{Y})$ d{\'e}fini par la
relation :

$$ \forall \sigma \in G_{\mathcal{U}},  \zeta_{\omega(\cal
  U)}(\sigma) = (\sigma-1)f.$$ C'est un cocycle de
$G_{\mathcal{U}}$ {\`a} valeurs dans $\omega(\mathcal{Y})$  dont la
restriction {\`a} $R_{u}$ est un morphisme de $C$-groupes alg{\'e}briques  injectif de $R_u$
  dans $\omega(\mathcal{Y})$.

\medskip
\begin{lemme}\label{lemme:norm}
L'image $W$ de $R_u$ sous $\zeta_{\omega(\cal U)}$ est un
sous-$G_{\mathcal{Y}}$-module de $\omega(\mathcal{Y})$.

\end{lemme}

\textbf{D{\'e}monstration}\\
Pour tout $\sigma_1  \in G_\mathcal{Y}$ et $\sigma_2 \in R_u$, on a
$$\zeta_{\omega(\cal U)}( \sigma_1 \sigma_2
{\sigma_1}^{-1})=\sigma_1(\zeta_{\omega(\cal U)}(\sigma_2)).$$

En effet, on  a : \begin{equation} \label{eqn:co2} \sigma_1
\zeta_{\omega(\cal U)}({\sigma_1}^{-1})= (1 - \sigma_1) f =
-\zeta_{\omega(\cal U)}(\sigma_1),\end{equation}\\
et
$$\zeta_{\omega(\cal U)}( \sigma_1 \sigma_2
{\sigma_1}^{-1})= \sigma_1(\zeta_{\omega(\cal
U)}(\sigma_2{\sigma_1}^{-1})) + \zeta_{\omega(\cal U)}(\sigma_1)
=\sigma_1(\sigma_2(\zeta_{\omega(\cal
U)}({\sigma_1}^{-1}))+\zeta_{\omega(\cal
U)}(\sigma_2))+\zeta_{\omega(\cal U)}(\sigma_1).$$ De
(\ref{eqn:co2}),  on  d{\'e}duit que :
$\sigma_1(\sigma_2(\zeta_{\omega(\cal
U)}({\sigma_1}^{-1})))=-\sigma_1 \sigma_2
{\sigma_1}^{-1}(\zeta_{\omega(\cal U)}(\sigma_1))$. Or $\sigma_1
\sigma_2 {\sigma_1}^{-1}$ est un {\'e}l{\'e}ment de
$Stab(\mathcal{Y})/Stab(\mathcal{U})$ et $\zeta_{\omega(\cal
U)}(\sigma_1)$ est dans $\omega(\mathcal{Y})$.\\
 On en d{\'e}duit
que $\sigma_1(\sigma_2(\zeta_{\omega(\cal U)}({\sigma_1}^{-1})))=
-\zeta_{\omega(\cal U)}(\sigma_1)$. Ainsi
$\sigma_1(\zeta_{\omega(\cal U)}(\sigma_2))=\zeta_{\omega(\cal U)}( \sigma_1 \sigma_2 {\sigma_1}^{-1})$ appartient bien {\`a} $W$.\\

\medskip

 \begin{prop}\label{prop:minimal} L' image sous $\omega$ de l'objet minimal de $\bold{V}(\mathcal{Y},\mathcal{U})$
co{\"\i}ncide avec  l'image sous  $\zeta_{\omega(\cal U)}$ de $R_u$.
\end{prop}

\textbf{D{\'e}monstration}\\
Notons
  $\mathcal{V}$ l'objet minimal de $\bold{V}(\mathcal{Y},\mathcal{U})$, et $V$ son image sous $\omega$. Alors, $G_\mathcal{U}$ agit sur $\omega(\mathcal{U}/\mathcal{V})$ {\`a} travers $G_\mathcal{Y}$ (puisque $\mathcal{U}/\mathcal{V}$ est une extension triviale de $ \bold{1}$ par un quotient de $\mathcal{Y}$). Donc la projection  de $f=s(1)$  sur $ \omega(\mathcal{U})/V$  est invariante sous $R_u$, et  l'orbite $\lbrace \sigma f -f ; \ \sigma \in R_u \rbrace$ est
 incluse dans $V$. Ainsi, $\zeta_{\omega(\cal U)} (R_u):= W  \subset V$.\\
R{\'e}ciproquement,  l'image  $W$ de $R_u$ sous $\zeta_{\omega(\cal U)}$  est, d'apr{\`e}s le  lemme
\ref{lemme:norm},
un sous-$G_{\mathcal{Y}}$-module de $\omega(\mathcal{Y})$ et est donc de la forme   $\omega(\mathcal{W})$ avec
$\mathcal{W} \subset \mathcal{Y}$.  Montrons
que $\mathcal{W} $ est un {\'e}l{\'e}ment de $\bold{V}$. Il est clair que  $Stab(\mathcal{Y})$
agit trivialement sur $\omega(\mathcal{U} / \mathcal{W})$ puisque pour tout $\sigma \in R_u = Stab({\mathcal{Y}})/Stab({\mathcal{U}})$ et $x \in \omega(\mathcal{U})$, $\sigma(x)-x \in
W = \omega(\mathcal{W})$.
Donc $\omega(\mathcal{U}/ \mathcal{W})$ est une repr{\'e}sentation
de $G_\mathcal{Y}$ et c'est une extension de $1_C$ par
$\omega(\mathcal{Y}/ \mathcal{W})$. Comme $G_\mathcal{Y}$ est
r{\'e}ductif,  cette extension se scinde comme extension dans la
cat{\'e}gorie $Rep_{G_{\mathcal{Y}}}$. D'apr{\`e}s la proposition \ref{prop:cocycle}, l'extension
$\mathcal{U} / \mathcal{W}$ est donc triviale dans
$Ext_{\bold{T}}(\bold{1}, \mathcal{Y} / \mathcal{W})$, et $\mathcal{W} \in \bold{V}$, d'o{\`u}   $V
\subset W$ par minimalit{\'e}. Ceci conclut la preuve du la proposition, et  du th{\'e}or{\`e}me \ref{theorem:Ber2}.

\medskip
\textbf{D{\'e}monstration  du corollaire \ref{coro:inde}} \\
Pour  toute extension $\mathcal{U}$ de $\bold{1}$ par $
\mathcal{Y}$, et tout $\alpha \in \Delta$, on d{\'e}signe par
$\alpha_* \mathcal{U}$   le \textit{pushout } de $\mathcal{U}$ par
$\alpha$  ; c'est ainsi que l'on d{\'e}finit la structure de
$\Delta$-module de $Ext^1_{\bold{T}}(\bold{1}, \mathcal{Y})$.

 L'extension $\mathcal{E}_1 \oplus ...\oplus \mathcal{E}_n$ de $\bold{1}^n$ par $\mathcal{Y}^n$ et son \textit{pull-back}  $\mathcal{E} \in Ext^1_{\bold{T}}(\bold{1}, \mathcal{Y}^n)$ par l'application diagonale de $\bold{1}$ dans $\bold{1}^n$ engendrent dans $\bold T$ la m{\^e}me sous-cat{\'e}gorie tannakienne, et admettent donc le m{\^e}me groupe de Galois $G_{\mathcal{E}_1 \oplus... \oplus
\mathcal{E}_n} = G_{\mathcal{E}}$.
Supposons que le radical unipotent de ce groupe  ne remplisse  pas
$\omega(\mathcal{Y}^n) = \omega(\mathcal{Y})^n$
  tout entier.\\
D'apr{\`e}s le th{\'e}or{\`e}me \ref{theorem:Ber2}, il existe un sous-objet
$\mathcal{V}$ non trivial de $\mathcal{Y}^n$ tel que le quotient de
$\mathcal{E}$ par $\mathcal{V}$ soit
  une  extension triviale de $\bold{1}$ par $\mathcal{Y}^n /
\mathcal{V}$. L'objet $\mathcal{Y}^n$ {\'e}tant compl{\`e}tement
r{\'e}ductible, il existe un {\'e}l{\'e}ment $\phi \in Hom(
\mathcal{Y}^n, \mathcal{Y})$ non nul dont le noyau $\mathcal{H}$
contienne $\mathcal{V}$ :

$$\xymatrix{
  {\mathcal{Y}}^n \ar[d] \ar[dr]^{\phi}  \\
{\mathcal{Y}}^n / \mathcal{V} \ar[r] & {\mathcal{Y}}^n / \mathcal{H}
\simeq  \mathcal{Y} .} $$ Ecrivons $\phi(X_1,...,X_n)= \alpha_1
X_1+...+\alpha_n X_n$, avec
  $\alpha_i \in End(\mathcal{Y})$. Alors $\phi_{*}(\mathcal{E})={\alpha_1}_* \mathcal{E}_1
+{\alpha_2}_* \mathcal{E}_2 +...+ {\alpha_n}_* \mathcal{E}_n$
est  un quotient de $\mathcal{E}/ \mathcal{V}$ donc  une extension triviale de $\bold{1}$ par $\mathcal{Y}$.\\
Ainsi l'extension $\alpha_1\mathcal{E}_1+...+\alpha_n
\mathcal{E}_n \in Ext^1(\bold{1},\mathcal{Y})$ est
triviale. Or  les extensions $\mathcal{E}_1,...,
\mathcal{E}_n$ sont $\Delta$-lin{\'e}airement ind{\'e}pendantes dans
$Ext^1_{\bold{T}}(\bold{1},\mathcal{Y}) $, d'o{\`u} la contradiction
attendue.

\subsection{Lien avec la transcendance}
\medskip
Pour pr{\'e}parer aux notations de la partie suivante, on d{\'e}signe ici par $C'$ le corps not{\'e} $C$ pr{\'e}c{\'e}demment. Rappelons qu'on ne le suppose pas alg{\'e}briquement clos. Comme dans l'introduction, on consid{\`e}re un corps  aux diff{\'e}rences, not{\'e} cette fois
$(K',\sigma)$, de corps des $\sigma$-constantes $C_{K'} := C'$. On d{\'e}signe par $Diff(K', \sigma)$ la cat{\'e}gorie des $K'[\sigma,\sigma^{-1}]$-modules  de dimension finie sur $K'$. C'est une cat{\'e}gorie tannakienne $C'$-lin{\'e}aire, qu'on suppose neutre dans l'{\'e}nonc{\'e} qui suit.\\

\begin{prop}\label{prop:degtr}
Soit $(K',\sigma)$ un corps aux diff{\'e}rences de corps des
constantes $C_{K'}=C'$
de caract{\'e}ristique nulle. \\
On suppose donn{\'e}e une extension $F$ de
corps aux diff{\'e}rences de $K'$, de corps des constantes
$C_F=C'$ v{\'e}rifiant la propri{\'e}t{\'e} suivante :  tout
syst{\`e}me aux $\sigma$-diff{\'e}rences   {\`a} coefficients dans $K'$ admet une matrice
fondamentale de solutions  {\`a} coefficients dans $F$.\\
 On note $\omega$ le foncteur fibre de $Diff(K',
\sigma)$ {\`a} valeurs dans $Vect_{C'}$ correspondant.

Soient $\mathcal{M} \in Diff(K', \sigma)$ et $U$ une matrice
fondamentale de solutions {\`a} coefficients dans $F$. On note
$G = Gal(\mathcal{M})$
 le groupe de Galois de $\mathcal{M}$ relativement {\`a}
 $\omega$.\\
 Alors
$$degtr(K'(U)/ K')= dim_{C'}G.$$
\end{prop}

\textbf{D{\'e}monstration}\\
La preuve qui suit {\'e}tend  la d{\'e}monstration de \cite{Ka} p.~38 au cas d'un corps de
constantes non alg{\'e}briquement clos.

Notons $<\mathcal{M}>$ la cat{\'e}gorie tannakienne engendr{\'e}e
par $\mathcal{M}$ dans $Diff(K',\sigma)$, de sorte que $G$ est le
$C'$-groupe alg{\'e}brique $Aut^{\otimes}(\omega|<\mathcal{M}>)$.
Soient $M \subset \mathcal{M}\otimes_{K'} F$ le $C'$-espace
vectoriel $\omega(\mathcal{M})$, et  $v$ un {\'e}l{\'e}ment de
$M(C')$. Alors,  $M$ est une $C'$-repr{\'e}sentation de $G$, et
  l'image sch{\'e}matique $G(v)$ du morphisme de $G$ dans $M$ attach{\'e} {\`a}
  $v$ est un $C'$-sous-sch{\'e}ma de $M$ (voir
\cite{Sp} 1.9.1 (iv)).  L'annulateur de $v \otimes 1$ dans $\mathcal{M}^*$ ($K'$-dual de
$\mathcal{M}$) est un sous-objet $\mathcal{W}$ de $\mathcal{M}^*$
dans la cat{\'e}gorie $<\mathcal{M}>$, et $W := \omega(\mathcal{W})$
est un sous-$G$-module du $C'$-dual $M^*$ de $M$. Comme $W \subset
\mathcal{W}\otimes_{K'} F$ annule $v$ et est $G$-stable, il doit
annuler $G(v)$. Ainsi, $W$ est contenu dans l'annulateur $S$ de
$G(v)$ dans $V^*$. Comme $S$ est stable sous $G$, il lui correspond
par {\'e}quivalence de cat{\'e}gories un sous-objet  $\mathcal{S}$
de $\mathcal{M}^*$ dans la cat{\'e}gorie $<\mathcal{M}>$. Comme $S$
annule $v$, il en est de m{\^e}me de $\mathcal{S} \subset
S\otimes_{C'} F$, et on a $\mathcal{S} \subset \mathcal{W}$. Mais $W
\subset S$, donc $\mathcal{S}=\mathcal{W}$. Ainsi, la dimension sur
$C'$ de l'espace des formes lin{\'e}aires sur $V$ annulant $G(v)$
est {\'e}gale {\`a}  la dimension sur $K'$ de l'espace des
formes lin{\'e}aires sur $\mathcal{M}$ qui annulent $v \otimes 1$.\\
En appliquant ce r{\'e}sultat au vecteur $\oplus_{i\leq
  n}Symm^i(v)$ et {\`a} l'objet $\oplus_{i\leq
  n}Symm^i(\mathcal{M})$ de \\
  $<\mathcal{M}>$, on obtient que le
degr{\'e} de transcendance de $K'(v)$ sur $K'$ est {\'e}gale {\`a} la dimension
sur $C'$  de $G(v)$.\\
Enfin, appliquons ce dernier r{\'e}sultat  {\`a} l'objet
$\mathcal{M}' := \underline{Hom}(M\otimes K', \mathcal{M})$ de
$<\mathcal{M}>$. Sa fibre $\omega(\mathcal{M}')$ s'identifie {\`a}
$End(M)$ et $G$ y agit par translation {\`a} gauche. La matrice
fondamentale de solutions $U$ de $\mathcal{M}$ correspond {\`a} un
{\'e}l{\'e}ment de $\omega(\mathcal{M}')$ dans $Gl(M)({C'})$. Comme
$G$ agit librement sur $Gl(M)/_{C'}$, on en d{\'e}duit que
$$degtr(K'(U)/ K')= dim_{C'}G.$$

\section{Changement de base }\label{sec:extscal}

Dans cette troisi{\`e}me partie, $C$ d{\'e}signe un corps
\textit{alg{\'e}briquement clos } de caract{\'e}ristique nulle, et $C'$ une extension \textit{quelconque} de $C$. En particulier,  \text\it {on ne suppose  pas que le corps $C'$ soit alg{\'e}briquement clos}.\\

 Soient $(K,\sigma)$ un corps aux diff{\'e}rences de corps des constantes $C$, et $(K',\sigma)$  une extension  aux diff{\'e}rences de $K$,   de corps des constantes
$C'$.
On d{\'e}signe par $ \mathcal{D}_K=K[\sigma, {\sigma}^{-1}]$ (resp. $ \mathcal{D}_{K'}=K'[\sigma,
{\sigma}^{-1}]$)  l'anneau des op{\'e}rateurs polynomiaux en
$\sigma$ et ${\sigma}^{-1}$ {\`a} coefficients dans $K$ (resp.
$K'$) et par $Diff(K, \sigma)$ (resp. $Diff(K', \sigma)$)  la
cat{\'e}gorie  des $\mathcal{D}_K$ (resp. $\mathcal{D}_{K'}$) modules de
dimension finie sur $K$ (resp. $K'$).\\
 On sait (voir \cite{Dlct} 6.20 ou \cite{VPS2}) (resp. on suppose) qu'il existe un foncteur fibre $\omega$ (resp. $\omega'$)
de la cat{\'e}gorie tannakienne $Diff(K, \sigma)$ (resp. $Diff(K',
\sigma)$)  vers la cat{\'e}gorie $Vect_C$ (resp. $Vect_{C'}$); en
d'autres termes,  que chacune de ces cat{\'e}gories tannakiennes est
neutre. Les   principaux cas qui nous int{\'e}resseront sont
d{\'e}velopp{\'e}s dans la partie  $4$ ($q$-diff{\'e}rences) et dans
la partie  $5$ ($\tau$-diff{\'e}rences), o{\`u} l'on verra que les
hypoth{\`e}ses suivantes sont satisfaites.

\begin{hyp}\label{hyp:ext}
\begin{enumerate}
\item  On suppose qu'il existe un morphisme de groupes injectif $i$
  d'un sous-groupe $\Gamma$  du
  groupe $G_C=Aut(C'/C)$ (groupe des
automorphismes de $C'$ sur  $C$  laissant $C$ invariant) dans
$G_K=Aut(K'/K)$ et on fixe un
  tel $i$. On se permet d'identifier   $\Gamma$ {\`a} son image $i(\Gamma)$.
\item $i$ fournit un  prolongement   de l'action naturelle de $\Gamma$ sur $C'$ {\`a}
  $K'$ et on suppose que cette action commute {\`a} l'action de $\sigma$
  sur $K'$ ; en particulier, $C'$ (resp. $C'^*$) est un
  sous-$\Gamma$-module de $K'$ (resp. $K'^*$).
\item On suppose que le corps des invariants de $K'$ sous l'action de $\Gamma$ est le
  corps $K$ lui-m{\^e}me: $K'^{\Gamma}= K$; en particulier  $C'^{\Gamma}=C$.
 \item On suppose enfin  que pour tout entier $n\geq 1$, l'application naturelle $H^1(\Gamma,Gl_n(C')) \rightarrow
 H^1(\Gamma, Gl_n(K')) $ est injective.

\end{enumerate}
\end{hyp}

Le r{\'e}sultat principal de cette partie $3$ sera d{\'e}montr{\'e}  au  paragraphe
$3.2$. Il {\'e}nonce :
\begin{theorem}\label{theorem:injext}
On suppose  les hypoth{\`e}ses \ref{hyp:ext} satisfaites.
Soit  $\mathcal{Y}$ un  objet de $Diff(K,
\sigma)$. L'application naturelle  de $Ext^1_{Diff(K, \sigma)}(\bold{1}, \mathcal{Y})$ dans
$Ext^1_{Diff(K', \sigma)}(\bold{1}, \mathcal{Y}\otimes K')$ est un morphisme de groupe injectif.
\end{theorem}

\subsection{Extension  du corps des constantes}
On reprend les notations du d{\'e}but de la partie  $3$ et on suppose
d{\'e}sormais que  \textit{les corps $K$ et $K'$ satisfont les hypoth{\`e}ses \ref{hyp:ext}}.

\subsubsection{Stabilit{\'e} des objets compl{\'e}tement r{\'e}ductibles}

\begin{lemme}\label{lemme:transir} Soit $\mathcal{M}$ un  $\mathcal{D}_K$-module irr{\'e}ductible. Alors
le $\mathcal{D}_{K'}$-module   $\mathcal{M} \otimes K'$ est compl{\`e}tement r{\'e}ductible.\\
\end{lemme}

\textbf{D{\'e}monstration}\\
Notons tout d'abord que si $\mathcal{N}$ est un $K$-sous-espace vectoriel de $ \mathcal{M}$ tel que $\mathcal{N} \otimes K'$ soit stable sous $\sigma$, alors $\mathcal{N}$ l'est aussi. En effet, $\sigma(\mathcal{N}\otimes 1) \subset (\mathcal{M}\otimes1) \cap (\mathcal{N}\otimes K') = \mathcal{N}\otimes 1 $.\\
Avec les  conventions de l'hypoth{\`e}se \ref{hyp:ext}.1,  consid{\'e}rons l'action semilin{\'e}aire de $\Gamma$ sur $\mathcal{M}
\otimes K'$ d{\'e}finie par $\tau \rightarrow id \otimes \tau$. Soient
$\mathcal{C}$ une composante isotypique maximale de
$\mathcal{M}\otimes K'$ dans ${Diff(K,
  \sigma)}$ et $\tau$ un {\'e}l{\'e}ment de $ \Gamma$. Alors $\tau(\mathcal{C})$ est {\`a} nouveau  une composante
  isotypique maximale de   $\mathcal{M} \otimes K'$ car l'action de $\Gamma$ et celle de
  $\sigma$ commutent. L'ensemble des composantes
  isotypiques maximales  de $\mathcal{M} \otimes K'$ {\'e}tant de cardinal fini, l'orbite de
  $\mathcal{C}$ sous $\Gamma$ est finie. Notons $(\tau_i(\mathcal{C}))_i$ les {\'e}l{\'e}ments
  de cette orbite; comme ce sont des composantes isotypiques maximales, elles sont lin{\'e}airement ind{\'e}pendantes. Posons alors $\mathcal{N}' = \bigoplus
  \tau_i(\mathcal{C})=\sum_{\tau \in \Gamma} \tau(\mathcal{C}) $.
  Le module
  $\mathcal{N}'$ est un sous-$\mathcal{D}_{K'}$-module de $\mathcal{M} \otimes K'$ stable sous l'action de
  $\Gamma$. D'apr{\`e}s l'hypoth{\`e}se \ref{hyp:ext}.3, il est  d{\'e}fini
    sur $K$, c'est-{\`a}-dire de la
  forme $\mathcal{N} \otimes K'$ o{\`u} $\mathcal{N} $est un $K$-sous espace vectoriel de $\mathcal{M}$, donc d'apr{\`e}s la remarque pr{\'e}liminaire,  un $\mathcal{D}_K$-sous-module de $\mathcal{M}$. Le
  $\mathcal{D}_K$-module  $\mathcal{M}$ {\'e}tant  irr{\'e}ductible, $\mathcal{N}=\mathcal{M}$. Donc
  $\mathcal{M}\otimes K' =\bigoplus \tau_i(\mathcal{C})$, somme directe de
  composantes isotypiques, est compl{\'e}tement r{\'e}ductible.\\

\begin{coro} \label{coro:transcomp} Soit $\mathcal{M}$ un  $\mathcal{ D}_K$-module compl{\'e}tement r{\'e}ductible. Alors
  $\mathcal{M} \otimes K'$ est un $\mathcal{D}_{K'}$-module compl{\`e}tement r{\'e}ductible.\\
\end{coro}

\textbf{D{\'e}monstration}\\
$\mathcal{M}$ est somme directe de  $\mathcal{D}_K$-modules
irr{\'e}ductibles $\mathcal{M}_i$. Donc $\mathcal{M}\otimes
  K' = \bigoplus \mathcal{M}_i \otimes K'$. D'apr{\`e}s le lemme \ref{lemme:transir}, $\mathcal{M}_i \otimes
  K'$ est compl{\'e}tement r{\'e}ductible, donc somme directe de  $\mathcal{D}_{K'}$-modules
  irr{\'e}ductibles $\mathcal{M}_i^j$. Donc $\mathcal{M}\otimes K' =\bigoplus \mathcal{M}_i^j$ est bien  somme
  directe d'objets irr{\'e}ductibles  de $Diff(K', \sigma)$.\\

\subsubsection{Compatibilit{\'e} des $Hom$}

\begin{definition} Soit $\mathcal{M}$ un objet de $Diff(K,
  \sigma)$. On dit que $\mathcal{M}$ est trivial sur $K$ s'il est
  isomorphe {\`a} une somme directe de copies de $\bold{1}$. Tout objet  $\mathcal{M}$ admet un plus grand sous-objet $\mathcal{N}$ trivial sur $K$,  et l'on a : $\mathcal{N} \simeq Hom(\bold{1},\mathcal{M}) \otimes_C K$.
\end{definition}

On montre dans ce paragraphe que sous les hypoth{\`e}ses
\ref{hyp:ext}, le plus grand sous-objet de
$\mathcal{M}'=\mathcal{M}\otimes_K K'$ trivial sur $K'$
(c'est-{\`a}-dire au sens de $Diff(K',
  \sigma)$)  provient par extension des scalaires de $K$ {\`a}
$K'$ du plus grand sous-objet de $\mathcal{M}$ trivial sur $K$ ou,
de fa{\c c}on {\'e}quivalente, que $Hom(\bold{1},\mathcal{M}) \otimes_C C'=
Hom(\bold{1},\mathcal{M}')$. De
fa{\c c}on g{\'e}n{\'e}rale, on peut {\'e}noncer:\\

\begin{lemme}\label{lemme:comphom}
On suppose que les hypoth{\`e}ses
\ref{hyp:ext} sont satisfaites. Soient $\mathcal{A}$ et $\mathcal{B}$ deux objets de
$Diff(K,\sigma)$. Alors:
$$ Hom_{Diff(K,
  \sigma)}(\mathcal{A},\mathcal{B}) \otimes_C  C' \simeq
Hom_{Diff(K',
  \sigma)}(\mathcal{A}\otimes_K K', \mathcal{B} \otimes_K K')$$
\end{lemme}

\textbf{D{\'e}monstration}\\
Posant $\mathcal{V}=\underline{Hom}(\mathcal{A},\mathcal{B})$, on se ram{\`e}ne {\`a} prouver que l'application naturelle $\iota$ de
\\ $ Hom_{Diff(K,
  \sigma)}(\bold{1},\mathcal{V}) \otimes_C  C'$ dans
$Hom_{Diff(K',
  \sigma)}(\bold{1}\otimes_K K', \mathcal{V} \otimes_K K')$, qui est une injection $C'$-lin{\'e}aire, est {\'e}galement surjective.\\
Notons  $\phi_{\sigma}$ l'action  de $\sigma $ sur
$\mathcal{V}$ et posons
$$N' = Hom_{Diff(K',
  \sigma)}(\bold{1}\otimes_K K', \mathcal{V}
\otimes_K K') = \lbrace g \in \mathcal{V}\otimes K' \ \mbox{tels que}
\ \phi_{\sigma}(g)=g \rbrace,$$ de sorte que $\mathcal{N}'=N'\otimes K'$ le plus grand sous-objet de $\mathcal{V}
\otimes_K K'$ trivial sur $K'$. Comme $\Gamma$ et $\sigma$ commutent, $\mathcal{N}'$ est de plus un
sous-module de $\mathcal{V}\otimes K'$
stable sous l'action de $\Gamma$.    D'apr{\`e}s l'hypoth{\`e}se
\ref{hyp:ext}.3,   $\mathcal{N}'$ est  donc d{\'e}fini
    sur $K$, c'est-{\`a}-dire de la
  forme $\mathcal{N} \otimes_K K'$, o{\`u} $\mathcal{N}$ un sous-$\mathcal{D}_K$-module de
    $\mathcal{V}$. On a :
$${\mathcal{N}'}^{\phi_{\sigma}} = N' = ({\mathcal{N} \otimes_K K'})^{\phi_{\sigma}}$$
Soient $n$ la dimension de $\mathcal{N}$ sur $K$, et $\sigma (Y) =R Y \ (*)$ une {\'e}quation aux diff{\'e}rences repr{\'e}sentant   $\mathcal{N}$ avec
$R \in Gl_n(K)$. Par d{\'e}finition de $N'$, il existe une matrice
fondamentale $U \in Gl_n(K')$ de solutions de $(*)$.\\
 Pour tout $\tau \in \Gamma$ , $\tau(U)$ est une solution de $(*)$.
 Ainsi :
 $$ \forall \tau \in \Gamma,  \ \exists R_{\tau} \in Gl_n(C')  \ \mbox{tel
 que} \  \tau(U)= U R_{\tau}.$$
Il est clair que $ \tau \mapsto R_{\tau}$ est un cocycle de $\Gamma$ {\`a}
valeur dans $Gl_n(C')$, trivial dans $H^1(\Gamma, Gl_n(K'))$.\\

D'apr{\`e}s l'hypoth{\`e}ses \ref{hyp:ext}.4: $H^1(\Gamma,
Gl_n(C')) \hookrightarrow H^1(\Gamma, Gl_n(K'))$,  il existe donc
un {\'e}l{\'e}ment $\beta$ de  $Gl_n(C')$ tel que pour tout  {\'e}l{\'e}ment $\tau$ de  $\Gamma$,
$B_{\tau}={\beta}^{-1} \tau(\beta)$. Alors:
$$ \forall \tau \in \Gamma, \tau(U{\beta}^{-1})= U{\beta}^{-1}.$$
Donc $U{\beta}^{-1}$ est un {\'e}l{\'e}ment de $Gl_n(K)$ solution de
$(*)$. En d'autres termes, il existe un   $C$-sous-espace vectoriel $N$ de dimension $n$ de
$Hom_{Diff(K,
  \sigma)}(\bold{1},\mathcal{V})$ tel que $N' = N \otimes_C C'$. Ainsi,   $\mathcal{N} =  N \otimes_C K$ est un sous-module de $\mathcal{V}$ trivial sur $K$. Comme $\iota$ est injective et que $N' = Hom_{Diff(K',
  \sigma)}(\bold{1}\otimes_K K', \mathcal{V}
\otimes_K K')$ est lui-m{\^e}me de dimension $n$ sur $C'$, on a n{\'e}cessairement $N =  Hom_{Diff(K,
  \sigma)}(\bold{1},\mathcal{V})$. La relation  $ N\otimes_C C'  = N'$ exprime alors  la surjectivit{\'e} de  l'application $\iota$.\\

\subsection{Injection des groupes d'extensions}

\subsubsection{D{\'e}monstration du th{\'e}or{\`e}me \ref{theorem:injext} }

 Soit  $\xymatrix{ 0 \ar[r] & \mathcal{X} \ar[r]^{i} & \mathcal{U}
    \ar[r]^{p} & \bold{1} \ar[r] & 0}$
une extension de $\bold{1}$ par $\mathcal{X}$ dans la cat{\'e}gorie $Diff(K,\sigma)$.\\
On note $\phi_{\sigma}$ l'action de $\sigma$ sur $\mathcal{U}$.\\
Dire que l'extension $\mathcal{U}$  se trivialise dans $K'$ (c'est-{\`a}-dire qu'elle appartient au noyau de l'application naturelle  de $Ext^1_{Diff(K, \sigma)}(\bold{1}, \mathcal{Y})$ dans
$Ext^1_{Diff(K', \sigma)}(\bold{1}, \mathcal{Y}\otimes K')$)
signifie qu'il existe un {\'e}l{\'e}ment $f \in \mathcal{U}\otimes
K'$ tel que
\begin{equation} \label{eqn:1} \phi_{\sigma}(f)=f \ \mbox{et} \ p(f) \neq 0
.\end{equation} Posons  $N'= \lbrace g \in \mathcal{U}\otimes K' \
    \mbox{tels que} \ \phi_{\sigma}(g)=g \rbrace$, de sorte que  $\mathcal{N}'=N'\otimes K'$ est
    le plus grand  sous-objet  de $\mathcal{U}\otimes K'$ trivial sur $K'$.
   D'apr{\`e}s le lemme \ref{lemme:comphom},  $\mathcal{N}'
    =\mathcal{N}\otimes K'$, o{\`u} $\mathcal{N}$ est le plus grand
    sous-objet de $\mathcal{U}$ trivial sur $K$.\\

    Ainsi, l'existence d'un {\'e}l{\'e}ment $f \in  \mathcal{N}'$ tel
    que $ p(f) \neq 0$ et $ \phi_{\sigma}(f)=f$ assure l'existence d'un {\'e}l{\'e}ment $g$ de $\mathcal{N} \subset \mathcal{U}$ tel
    que $p(g) \neq 0$ et $\phi_{\sigma}(g)=g$. L'extension $\mathcal{U}$ est donc alors triviale dans $Ext^1_{Diff(K, \sigma)}(\bold{1}, \mathcal{Y})$.\\

\subsubsection{Illustrations}

\medskip
Joint au th{\'e}or{\`e}me \ref{theorem:injext}, le corollaire \ref{coro:inde}, appliqu{\'e} {\`a} la cat{\'e}gorie $Diff(K', \sigma)$,  entra{\^\i}ne imm{\'e}diatement:\\

\begin{coro}\label{coro:indextg}
Soient  $\mathcal{Y}$ un objet
  compl{\'e}tement r{\'e}ductible   de $Diff(K,\sigma)$,  $\Delta$ l'anneau $End(\mathcal{Y})$  et
 $\mathcal{E}_1,...,\mathcal{E}_n$ des extensions de $\bold{1}$ par $\mathcal{Y}$, qu'on suppose  $\Delta$-lin{\'e}airement ind{\'e}pendantes dans $Ext^1_{Diff(K, \sigma)} (\bold{1},\mathcal{Y})$. Alors le
radical unipotent de $G_{(\mathcal{E}_1 \oplus... \oplus
  \mathcal{E}_n)\otimes K'}$ est isomorphe {\`a}
  $(\omega(\mathcal{Y})\otimes C')^n$.
\end{coro}

 \medskip
Voici une illustration concr{\`e}te de la fa\c {c}on de v{\'e}rifier les hypoth{\`e}ses  de ce corollaire. Rappelons que pour $ \mathcal{Y} = \bold{1}$, on a   $End_{Diff(K, \sigma)}(\bold{1}) = C$ et $End_{Diff(K', \sigma)}(\bold{1}) = C'$.

\medskip
\begin{prop}\label{prop:indext1}
 Soient $b_1, ..., b_n$ des {\'e}l{\'e}ments de $K$, et pour tout $i = 1, ..., n$,  $\mathcal{E}_i$  l'extension de $\bold{1}$
 par $\bold{1}$ dans la cat{\'e}gorie $Diff(K,\sigma)$, de repr{\'e}sentation matricielle donn{\'e}e par $\left( \begin{array}{cc}
 1 & b_i \\
 0 & 1 \end{array} \right)$. Les propri{\'e}t{\'e}s suivantes sont {\'e}quivalentes:

 i)  les
 extensions $\mathcal{E}_i$
 sont $C$-lin{\'e}airement ind{\'e}pendantes dans
 $Ext^1_{Diff(K, \sigma)}(\bold{1},\bold{1})$;

 ii) les {\'e}l{\'e}ments  $b_i$ sont $C$-lin{\'e}airement ind{\'e}pendants modulo $(\sigma_q
 -id)(K)$;

iii) les
 extensions $\mathcal{E}_i\otimes K'$
 sont $C'$-lin{\'e}airement ind{\'e}pendantes dans
 $Ext^1_{Diff(K', \sigma)}(\bold{1},\bold{1})$;

 iv) le groupe de Galois $G_{(\mathcal{E}_1 \oplus... \oplus
  \mathcal{E}_n)\otimes K'}$ est isomorphe {\`a} $\bold{G}_a^n/_{C'}$.

  \noindent
  Plus g{\'e}n{\'e}ralement, soit $\delta$ la dimension du $C$-sous espace vectoriel engendr{\'e} par $b_1, ..., b_n$ dans $K/(\sigma- 1)(K)$. Alors,  $G_{(\mathcal{E}_1 \oplus... \oplus
  \mathcal{E}_n)\otimes K'} \simeq \bold{G}_a^\delta/_{C'}.$
 \end{prop}

 \textbf{D{\'e}monstration}\\
 L'{\'e}quivalence de (i) et (iii) r{\'e}sulte du th{\'e}or{\`e}me  \ref{theorem:injext}. Pour celle de (i) et (ii), noter que pour tout $n$-uplet $\alpha_1, ..., \alpha_n$ d'{\'e}l{\'e}ments de $C$, l'extension $  \alpha_1\mathcal{E}_1 + ... + \alpha_n\mathcal{E}_n$ admet pour repr{\'e}sentation matricielle $\left( \begin{array}{cc}
 1 & b \\
 0 & 1 \end{array} \right)$, o{\`u} $b = \alpha_1b_1 + ... + \alpha_n b_n$. Une telle extension correspond {\`a} l'{\'e}quation aux diff{\'e}rences non homog{\`e}ne $\sigma y - y = b$, et est triviale dans  $Ext^1_{Diff(K, \sigma)}(\bold{1},\bold{1})$ si et seulement si $b$ appartient {\`a} l'image de $K$ sous l'op{\'e}rateur $\sigma -1$. \\
 Puisque $\mathcal{Y} = \bold{1}$ est ici trivial, le groupe de Galois $G_{(\mathcal{E}_1 \oplus... \oplus
  \mathcal{E}_n)\otimes K'} $ co{\"\i}ncide avec son radical unipotent, et l'implication $(i)$ (ou $(iii) ) \Rightarrow (iv)$ d{\'e}coule du corollaire \ref{coro:indextg}. La r{\'e}ciproque se d{\'e}duit de l'{\'e}nonc{\'e} plus pr{\'e}cis fourni par la proposition \ref {prop:minimal}.\\
  Soit enfin $\{b_{i_1}, ..., b_{i_\delta}\}$  un syst{\`e}me maximal extrait de  $\{b_1, ..., b_n\} $ et form{\'e} d'{\'e}l{\'e}ments $C$-lin{\'e}airement ind{\'e}pendants modulo $(\sigma-1)(K)$. Alors, $\mathcal{E}_{i_1} \oplus   ... \oplus \mathcal{E}_{i_\delta}$ engendre la sous-cat{\'e}gorie $<\mathcal{E}_1 \oplus ... \oplus \mathcal{E}_n>$ de $Diff(K, \sigma_q)$ , et la derni{\`e}re assertion d{\'e}coule des {\'e}quivalences pr{\'e}c{\'e}dentes, appliqu{\'e}es {\`a} ces $\delta$ extensions.

\subsection{Comparaison des groupes de Galois aux diff{\'e}rences}

Ce paragraphe ne sera pas utilis{\'e} par la suite. On y donne une nouvelle d{\'e}monstration du corollaire \ref{coro:indextg}, par le biais d'un {\'e}nonc{\'e} plus intrins{\`e}que (th{\'e}or{\`e}me \ref{theorem:extscal}),  qui permet de limiter les calculs des radicaux unipotents $R_u$ des groupes de
Galois {\`a} ceux de la cat{\'e}gorie $\bold{T}$.  Ce th{\'e}or{\`e}me \ref{theorem:extscal} est, pour l'essentiel, un cas particulier du th{\'e}or{\`e}me de comparaison g{\'e}n{\'e}ral {\'e}tabli dans  \cite{CHS}, en vertu duquel, pour tout objet $\mathcal{M}$ de  $Diff(K,
\sigma)$ et tout  couple de foncteurs  fibres  $\omega, \omega'$ de $Diff(K,
\sigma), Diff(K', \sigma)$ sur $C, C'$ (notations du d{\'e}but du
paragraphe 3), on a, sous certaines conditions sur l'extension $K'/K$:
\begin{equation}\label{chs} Aut^{\otimes}(\omega|<\mathcal{M}>)\otimes_C \overline{C'}
 \simeq  Aut^{\otimes}( \omega'|< \mathcal{M}\otimes K' >) \otimes_{C'} \overline{C'}
 \end{equation}
o{\`u} $\overline{C'}$ d{\'e}signe une cl{\^o}ture alg{\'e}brique de $C'$.

On se place maintenant sous les hypoth{\`e}ses \ref{hyp:ext}.


\begin{theorem} \label{theorem:extscal}
On suppose que les hypoth{\`e}ses \ref{hyp:ext}   sont satisfaites.
Soient  $\mathcal{Y}$ un  objet compl{\'e}tement r{\'e}ductible de
$Diff(K, \sigma)$ et $\mathcal{U}$ une extension dans $Diff(K,
\sigma)$ de $\bold{1}$ par $\mathcal{Y}$. Alors
$$ R_u(Aut^{\otimes}( \omega{\^E}| < \mathcal{U} > )) \otimes C'
\simeq R_u( Aut^{\otimes}( \omega' | < \mathcal{U}\otimes K' >) ).$$
\end{theorem}
Gr{\^a}ce {\`a} cet {\'e}nonc{\'e},  le corollaire \ref{coro:indextg} devient une cons{\'e}quence imm{\'e}diate du corollaire \ref{coro:inde}, appliqu{\'e} cette fois {\`a} la cat{\'e}gorie $Diff(K, \sigma)$: en effet,  selon ce dernier, $R_u(G_{\mathcal{E}_1 \oplus... \oplus
\mathcal{E}_n}) \simeq \omega(\mathcal{Y})^n$; le th{\'e}or{\`e}me \ref{theorem:extscal} entra{\^\i}ne alors que $R_u(G_{(\mathcal{E}_1 \oplus... \oplus
\mathcal{E}_n)\otimes K'}) \simeq \omega(\mathcal{Y})^n\otimes C' .$

\begin{remark}
La principale hypoth{\`e}se sur $K'/K$ faite   dans  \cite{CHS} est que $K' = K(C')$; tout comme les hypoth{\`e}ses \ref{hyp:ext}, cette hypoth{\`e}se sera remplie dans les exemples des parties 4 et 5. On notera par ailleurs que la
  comparaison des groupes de Galois de \cite{CHS} impose de passer {\`a}
  la cl{\^o}ture alg{\'e}brique de $C'$. Dans le cas des groupes unipotents ab{\'e}liens, on peut toutefois, en
  utilisant \cite{Sp}  $12.3.8$, redescendre
l'isomorphisme (\ref{chs})  jusqu'{\`a} $C'$.
 \end{remark}

\textbf{D{\'e}monstration du th{\'e}or{\`e}me   \ref{theorem:extscal}}\\
D'apr{\`e}s le th{\'e}or{\`e}me \ref{theorem:Ber2}, le radical
unipotent du groupe de Galois de $\mathcal{U}$ (resp.
$\mathcal{U}\otimes K'$)  est {\'e}gal au  groupe vectoriel
$\omega(\mathcal{V})$ (resp. $\omega'(\mathcal{V}')$), o{\`u}
$\mathcal{V}$ (resp.$\mathcal{V}'$)  est le plus petit sous-objet de
$\mathcal{Y}$ dans  $Diff(K,\sigma)$ (resp. de $\mathcal{Y}\otimes
K'$ dans  $Diff(K',\sigma)$) poss{\'e}dant la propri{\'e}t{\'e}
suivante:  le quotient  de l'extension $\mathcal{U}$ (resp.
$\mathcal{U} \otimes K'$)  par $\mathcal{V}$ (resp.  $\mathcal{V}'$)
est une extension  scind{\'e}e. On en d{\'e}duit que  le quotient
de l'extension $\mathcal{U} \otimes K'$ par $\mathcal{V}\otimes K'$
dans la
cat{\'e}gorie $Diff(K',\sigma)$  est une extension  scind{\'e}e. Par cons{\'e}quent, $\mathcal{V}' \subset \mathcal{V}\otimes K'$.\\

R{\'e}ciproquement, consid{\'e}rons l'action semilin{\'e}aire de
$\Gamma$ ($\tau \rightarrow id \otimes \tau$) sur $\mathcal{U}
\otimes K'$. Cette action permute les sous
$\mathcal{D}_{K'}$-modules de $\mathcal{U} \otimes K'$. Soit $\tau
\in \Gamma$. Le quotient de l'extension  $\mathcal{U} \otimes K'$
par $\tau(\mathcal{V}')$ est une extension scind{\'e}e. Donc
$\mathcal{V}' \subset  \tau( \mathcal{V}')$ pour tout $\tau \in
\Gamma$ ; comme les dimensions sur $K'$ de ces deux modules sont
{\'e}gales, on a  $\tau(\mathcal{V}') =\mathcal{V}'$.
L'hypoth{\`e}se \ref{hyp:ext}.3 entra{\^\i}ne alors que
$\mathcal{V}'$ est d{\'e}fini sur $K$, et il existe un sous
$\mathcal{D}_{K}$-module $\mathcal{W}$ de $\mathcal{U}$ tel que
$\mathcal{V}'= \mathcal{W} \otimes K'$. En particulier, l'extension
$(\mathcal{U}\otimes K' )/ \mathcal{V'}$ de $\bold 1$ par
$(\mathcal{Y} \otimes K')/ \mathcal{V'}$ se d{\'e}duit par
changement de base de $K$ {\`a} $K'$ de l'extension
$\mathcal{U}/\mathcal{W}$ de $\bold 1$ par $\mathcal{Y} /
\mathcal{W}$, et on d{\'e}duit du th{\'e}or{\`e}me
\ref{theorem:injext} que cette derni{\`e}re est scind{\'e}e dans la
cat{\'e}gorie $Diff(K,\sigma)$. Ainsi $\mathcal{V} \subset
\mathcal{W}$.

On a donc $$\mathcal{V}'= \mathcal{W} \otimes K'=\mathcal{V}\otimes
K'.$$ D'apr{\`e}s le th{\'e}or{\`e}me \ref{theorem:Ber2}, $
R_u(Aut^{\otimes}( \omega{\^E}| < \mathcal{U} > )) \simeq
\omega(\mathcal{V})$ et $R_u( Aut^{\otimes}( \omega' | <
\mathcal{U}\otimes K' >) ) \simeq \omega'(\mathcal{V}') =
\omega'(\mathcal{V}\otimes K' )$. Comme $\omega(\mathcal{V}) \otimes
C'$ et $\omega'(\mathcal{V}\otimes K')$  ont m{\^e}me dimension sur
$C'$, cela termine la d{\'e}monstration du th{\'e}or{\`e}me
\ref{theorem:extscal}.

\section{Hypertranscendance des solutions d'{\'e}quations aux $q$-diff{\'e}rences}

Soit $q \in \mathbb{C}^*$ un nombre complexe de module diff{\'e}rent
de $1$. On d{\'e}signe par $K=\mathbb{C}(z)$ le corps des fractions rationelles
{\`a} coefficients complexes, par $F=\mathcal{M}er(\mathbb{C}^*)$ le corps
des fonctions m{\'e}romorphes sur $\mathbb{C}^*$ et par $\sigma_q$
l'automorphisme de $F$ qui {\`a} $f(z) \in F$ associe $f(qz)$. Comme dans l'introduction, on note $C_E$
 l'ensemble des fonctions de $F$ fix{\'e}es par
${\sigma}_q$, et $K_E=C_E(z)$ le compositum de $C_E$ et de $K$ dans
$F$. Les corps $K, F$ et $K_E$ sont des corps aux diff{\'e}rences
relativement {\`a} $\sigma_q$, admettant respectivement pour corps
des $\sigma_q$-constantes $\mathbb{C},  C_E$ et $C_E$, et $C_E$
s'identifie au corps des fonctions rationnelles sur la
courbe elliptique $E=\mathbb{C}^*/q^{\mathbb{Z}}$.\\

On d{\'e}signe par $ \mathcal{D}_K=K[\sigma_q, {\sigma_q}^{-1}]$ (resp. $ \mathcal{D}_{K_E}=K_E[\sigma_q,
{\sigma_q}^{-1}]$)  l'anneau des op{\'e}rateurs polynomiaux en
$\sigma_q$ et ${\sigma_q}^{-1}$ {\`a} coefficients dans $K$ (resp. dans
$K_E$), et par $Diff(K, \sigma_q)$ (resp. $Diff(K_E, \sigma_q)$)   la
cat{\'e}gorie  tannakienne $\mathbb{C}$-lin{\'e}aire (resp. $C_E$-lin{\'e}aire) des $\mathcal{D}_K$ (resp. $\mathcal{D}_{K_E}$) modules de
dimension finie sur $K$ (resp. $K_E$).\\

Dans \cite{VPS2}, M.F. Singer et M. van der Put construisent un
foncteur fibre $\omega$ de la cat{\'e}gorie $Diff(K, \sigma_q)$
{\`a} valeur dans la cat{\'e}gorie $Vect_{\mathbb{C}}$ des espaces
vectoriels de dimension finie sur $\mathbb{C}$. Ce foncteur fibre
attache {\`a} tout module aux $q$-diff{\'e}rences $\mathcal{M}$ sur
$K$ une $\mathbb{C}$-base de solutions symboliques et munit $ Diff(K,
\sigma_q)$ d'une structure de cat{\'e}gorie tannakienne  neutre sur $\mathbb{C}$.\\

 Dans \cite{Prag}, C. Praagman
d{\'e}montre que toute {\'e}quation aux $q$-diff{\'e}rences {\`a}
coefficients m{\'e}romorphes sur $\mathbb{C}^*$ admet une $C_E$-base
de solutions m{\'e}romorphes sur $\mathbb{C}^*$. Ceci permet de
construire un foncteur fibre $\omega_E$ de $Diff(K_E, \sigma_q)$
{\`a} valeur dans $Vect_{C_E}$, qui attache {\`a} tout module aux
  $q$-diff{\'e}rences sur $K_E$ le  $C_E$-espace vectoriel   de ses solutions dans $F=\mathcal{M}er(\mathbb{C}^*)$ et
  de munir $Diff(K_E,\sigma_q)$ d'une structure de cat{\'e}gorie
  tannakienne neutre sur le corps (non alg{\'e}briquement clos) $C_E$. Comme on l'a dit au d{\'e}but de l'article, c'est ce foncteur fibre que nous privil{\'e}gierons, car contrairement {\`a} celui de \cite{VPS2}, il conduit naturellement {\`a} des extensions de corps aux  diff{\'e}rences {\it diff{\'e}rentiels}.

 Plus pr{\'e}cis{\'e}ment, l'automorphisme $\sigma_q$ et la d{\'e}rivation $\partial = zd/dz$
munissent le corps $F =  \calM er(\mathbb{C}^*)$ et ses sous-corps
$K_E$ et $K$ de structures de corps aux diff{\'e}rences
diff{\'e}rentiels, puisque $\sigma_q \partial = \partial \sigma_q$
(en effet, pour tout $f \in F$, on a $\sigma_q
\partial f(z)= \sigma_q(z\frac{d}{dz}(f(z))=qz
(\frac{d}{dz}f)(qz)=z\frac{d}{dz}(f(qz))= \partial \sigma_q f(z)$).
Il en est de m{\^e}me du sous-corps de $F$ engendr{\'e} sur $K_E$
par toutes les d{\'e}riv{\'e}es $\partial^m Y, m \geq 0$ de la
solution $Y$ de l'{\'e}quation (\ref{eqn:der1q}) consid{\'e}r{\'e}e
dans l'introduction.  {\`A} l'exception de $K$, ces corps admettent
chacun $C_E$ comme corps de $\sigma_q$-constantes. L'accent portera
donc ici sur le  corps aux diff{\'e}rences diff{\'e}rentiel
$(K_E:=C_E(z),\sigma_q,\partial =z\frac{d}{dz})$, et sur la
cat{\'e}gorie  ${Diff(K_E, \sigma_q)}$. En particulier, les groupes
de Galois aux diff{\'e}rences $G=Gal_{\omega_E}$ calcul{\'e}s plus
bas sont relatifs au foncteur fibre
$\omega_E$ de Praagman, et  sont donc des groupes alg{\'e}briques sur $C_E$.\\

Nous sommes  maintenant en mesure de d{\'e}montrer le th{\'e}or{\`e}me
\ref{theorem:gen}. On d{\'e}montrera ce th{\'e}or{\`e}me dans le cadre particulier d'un objet de
 ${Diff(K_E, \sigma_q)}$ de rang $1$ au th{\'e}or{\`e}me \ref{theorem:hyper}
 et dans son cadre  g{\'e}n{\'e}ral au th{\'e}or{\`e}me \ref{theorem:gen2}.\\

\begin{Not}\label{Not:ext}
Si $\mathcal{M}$ est  un objet de $Diff(K, \sigma_q)$, on notera
$\mathcal{M}_E = \mathcal{M}\otimes_K K_E$ l'objet de $Diff(K_E,
\sigma_q)$ qu'on en d{\'e}duit par changement de base, et
$Gal_{\omega_E}(\cal M)$ le groupe  de Galois de $ \mathcal{M}_E$
relativement au foncteur fibre $\omega_E$. Enfin, on notera
$K_E(\mathcal{M})$ le corps de d{\'e}finition des {\'e}l{\'e}ments
de $ \omega_E(\mathcal{M}_E) \subset \mathcal{M}_E \otimes_{K_E}F$
relativement {\`a} la $K_E$-structure de $\mathcal{M}_E$,
c'est-{\`a}-dire le corps de d{\'e}finition, not{\'e} $K_E(U)$ {\`a}
la proposition \ref{prop:degtr}, d'une matrice fondamentale de
solutions $U$ de $\mathcal{M}_E$ dans $F$.
\end{Not}

\subsection{V{\'e}rification des hypoth{\`e}ses  \ref{hyp:ext}}\label{subsec:verhq}

On se propose de v{\'e}rifier ces hypoth{\`e}ses pour le couple de corps aux diff{\'e}rences form{\'e} de $K=\mathbb{C}(z)$ et de son extension  $K' = K_E = C_E(z)$. On note   $G_E = \Aut(C_E/\mathbb{C}) $ (resp.  $G_K=Aut(K_E/K)$) le groupe des automorphismes  de l'extension
$C_E/\mathbb{C}$ (resp. $K_E/K$), et $C_E(X)$ le corps des fractions
  rationnelles en une variable {\`a} coefficients dans  $C_E$.  Pour tout {\'e}l{\'e}ment $a$ de $E(\mathbb{C})$, la translation par $a$ d{\'e}finit un automorphisme de la courbe $E$; on note $\gamma_a \in G_E$ l'automorphisme correspondant de l'extension $C_E/\mathbb{C}$, ainsi que son prolongement canonique au corps $C_E(X)$, d{\'e}fini par son action sur les coefficients. L'ensemble $\Gamma = \lbrace \gamma_a, a \in E(\mathbb{C}) \rbrace$ forme un sous-groupe  de $G_E$, isomorphe {\`a} $E(\mathbb{C})$.

\begin{lemme}\label{lemme:galois} Avec  ces notations,
\begin{enumerate}

\item il existe un isomorphisme canonique $\phi$ de $ C_E(X)$ sur $K_E$, induisant un isomorphisme de   $\mathbb{C}(X)$ sur $K$;
\item soit   $a$ un  point de la courbe elliptique $E$. Par transport de structure par $\phi$, l'automorphisme $\gamma_a$ de $C_E(X)$ d{\'e}finit un automorphisme  $g_a$  de  l'extension $K_E/K$,  et l'application $i: \gamma_a \mapsto g_a$ est un homomorphisme injectif  de $\Gamma$ dans $G_K$;
\item   $ C_E^{\Gamma} =\mathbb{C}$, et $ K_E^{\Gamma} = K$;
\item l'application naturelle $H^1(\Gamma,Gl_n(C_E)) \rightarrow H^1(\Gamma,Gl_n( K_E))$ est injective;
\item l'action de $\Gamma$ sur $K_E$ commute  {\`a} celle  de $\sigma_q$.
\end{enumerate}
\end{lemme}

\textbf{D{\'e}monstration}\\
\begin{enumerate}

\item  Pour tout   $f(X) \in C_E[X]$, posons  $\phi(f) = f(z)$, vu comme
une fonction m{\'e}romorphe en $z \in  \mathbb{C}^*$. Alors , $\phi$ est
un morphisme de  $C_E[X]$ sur
  $K_E$. Par ailleurs, les sous-corps $C_E$ et $\mathbb{C}(z)$ de $F$
  sont lin{\'e}airement disjoints sur $\mathbb{C}$ (voir \cite{Eti}, p.~5,
  premier lemmme). Par cons{\'e}quent,
  toute relation de d{\'e}pendance dans  $\mathcal{M}er(\mathbb{C}^*)$ :
\begin{equation} \label{eqn:22} \sum c_i(z)k_i(z)=0, \forall z \in \mathbb{C}^*
\end{equation}
avec  $c_i \in C_E$ et  $k_i \in K$ entra{\^\i}ne:
\begin{equation} \label{eqn:3} \sum c_i(z)k_i(X)=0, \forall z \in
  \mathbb{C}^*. \\
\end{equation}
L'application $\phi$ est donc injective. Ainsi,  $\phi$ s'{\'e}tend au corps des fractions  $C_E(X)$, et son
image est $K_E$ tout entier. Remarquons que, par d{\'e}finition de
$\phi$,
 $\mathbb{C}(X)$ s'envoie isomorphiquement sur  $K$.\\

\item  Ceci d{\'e}coule de la d{\'e}finition de l' action de  $\Gamma$ sur
  $K_E$.  Comme   $\Gamma$ agit trivialement sur
  $\mathbb{C}(X)$, son action sur   $K$ est aussi triviale.\\

\item

Une fonction elliptique invariante par toutes les translations est constante, donc $C_E^{\Gamma} = \mathbb{C}$. Par cons{\'e}quent, $C_E(X)^{\Gamma} = \mathbb{C}(X)$, et   $K_E^{\Gamma}=K$.\\

\item
 Soit $a \mapsto c_a$
un cocycle de $\Gamma \simeq E(\mathbb{C})$ {\`a} valeurs dans $Gl_n(C_E)$ trivial dans
$H^1(\Gamma, Gl_n(K_E)) $. Alors, il existe   $$ A(z)=
(P_{i,j}(z)=\frac{\sum_{k=0}^{n_{i,j}}
c^k_{i,j}(z)z^k}{\sum_{k=0}^{n'_{i,j}} c'^k_{i,j}(z)z^k} )_{i,j}$$
avec $c^k_{i,j},d^k_{i,j}, c'^k_{i,j},d'^k_{i,j}  \in C_E$, et
$A(z)  \in Gl_n(K_E)$ tels que :
\begin{equation}\label{eqn:trivcq} \forall a \in E(\mathbb{C}) , c_a(z)
  = A(z)^{-1} A(az).\end{equation}
Posons
$$P_{i,j}(z,X)=\frac{\sum_{k=0}^{n_{i,j}}
c^k_{i,j}(z)X^k}{\sum_{k=0}^{n'_{i,j}}
c'^k_{i,j}(z)X^k}$$
et $A(z,X) =
(P_{i,j}(z,X))_{i,j} \in GL_n(C_E(X))$

L'application $\phi$ {\'e}tant injective,  la relation (\ref{eqn:trivcq})
implique que :
\begin{equation}\label{eqn:trivcq1} \forall a \in E(\mathbb{C}) ,  \  c_a(z)
  = A(z,X)^{-1} A_2(az,X) .\end{equation}
Les fractions rationnelles $P_{i,j}(z,X),
det(A(z,X))  \in C_E(X)$ n'ayant qu'un nombre fini
de z{\'e}ros et de p{\^o}les, il existe un nombre complexe   $t_0$ tel
que pour tout $i, j$, et tout $a$ dans  $E(\mathbb{C})$,
$$ P_{i,j}(az,t_0) , det(A(az,t_0)) \neq 0, \ \infty  \\
~\mbox{dans} \ C_E.$$
L'action de $g_a$  sur $K_E \simeq C_E(X)$
commutant avec la
  sp{\'e}cialisation en $t_0$,  on d{\'e}duit que
\begin{equation}\label{eqn:trivcq2} \forall a \in E(\mathbb{C}), c_a(z)
  =A(z,t_0)^{-1}A(az,t_0).\end{equation}

Cette derni{\`e}re relation exprime que le cocycle $c_a$ est un
cobord {\`a}
valeurs dans $Gl_n(C_E)$, ce qui termine la d{\'e}monstration.\\

\item   Soit  $ k$ un entier naturel et  $f(X) = c X^k$ avec  $c \in
  C_E$. Alors
$$ \tau ( \sigma (f)) =\tau (c q^k X^k)= \tau(c) q^k X^k
=\sigma(\tau (f))$$
pour tout  $\tau \in \Gamma$. Ainsi, l'action  de $\Gamma$
commute avec  $\sigma$ sur  $C_E[X]$. Ces deux actions commutent
donc sur  $C_E(X) \simeq K_E$.\\

\end{enumerate}

\textbf{Extension des scalaires et $q$-diff{\'e}rences}\\

Le couple $(K,K_E)$ v{\'e}rifie les hypoth{\`e}ses de la partie \ref{sec:extscal}. On
peut donc appliquer le th{\'e}or{\`e}me \ref{theorem:injext} au cas des
$q$-diff{\'e}rences. On obtient, en reprenant les notations  \ref{Not:ext} :

\begin{theorem}\label{theorem:qinjext}
Soit  $\mathcal{Y}$ un  objet de $Diff(K,
\sigma_q)$. L'application naturelle  de $Ext^1_{Diff(K, \sigma_q)}(\bold{1}, \mathcal{Y})$ dans
$Ext^1_{Diff(K_E, \sigma_q)}(\bold{1}, \mathcal{Y}_E)$ est un morphisme de groupe injectif.
\end{theorem}

De m{\^e}me la proposition \ref{prop:indext1}   entra{\^\i}ne, pour les groupes de Galois attach{\'e}s au foncteur fibre $\omega_E$ de ${Diff(K_E, \sigma_q)}$:

\begin{coro}\label{coro:qext}
Soient $b_1, ..., b_n$ des {\'e}l{\'e}ments de $K$, et pour tout $i = 1, ..., n$,  $\mathcal{E}_i$  l'extension de $\bold{1}$
 par $\bold{1}$ dans la cat{\'e}gorie $Diff(K,\sigma)$, de repr{\'e}sentation matricielle donn{\'e}e par $\left( \begin{array}{cc}
 1 & b_i \\
 0 & 1 \end{array} \right)$. Alors le
groupe de Galois $G_{\omega_E}(\mathcal{E}_1 \oplus... \oplus
  \mathcal{E}_n)$ est isomorphe {\`a}
$\bold{G}_a^\delta/_{C_E}$, o{\`u} $\delta$ d{\'e}signe la dimension du $\mathbb{C}$-espace vectoriel engendr{\'e} par $b_1, ..., b_n$ dans $K/(\sigma_q -1)(K)$. .
\end{coro}

\subsection{Equations du premier ordre}

\begin{definition}\label{def:stand}
Soit  $a \in K^*=\mathbb{C}(z)^*$. On dira que $a$ est  standard
(\cite{VPS2} chap.2) si, pour tout $c \in {\mathbb{C}}^*$ dans le
support du diviseur de $a$, et tout entier $m$ non nul, $q^m c$
n'appara{\^\i}t pas dans le diviseur de $a$.
\end{definition}
\begin{lemme}\label{lemme:qstand}
Soit $a$ un {\'e}l{\'e}ment de $K^*$.
\begin{enumerate}
\item Il existe un  couple $(g,  \overline{a}$), avec $g
\in K^*$ et  $\overline{a}$  standard, tel que  $a = \overline{a}
\frac{\sigma_q(g)}{g}$. Une telle d{\'e}composition est dite forme
standard de $a$.
\item Soient $f$ (resp. $\overline{f}$) une solution
m{\'e}romorphe sur $\mathbb{C}^*$ de $\sigma_q(f)=af$ (resp.
$\sigma_q(\overline{f})=\overline{a}\overline{f}$). Pour tout $n
\geq 0$,  les corps $K(\partial^i(f);0 \leq  i \leq n)$ et
$K(\partial^i(\overline{f}); 0 \leq  i \leq n)$ co{\"\i}ncident.
\end{enumerate}
\end{lemme}

\textbf{D{\'e}monstration}\\
Pour le premier point, voir \cite{VPS2} p~29 lemme 2.1 et 2.2, qui le traite dans le cadre
des $\tau$-diff{\'e}rences. La d{\'e}monstration est la m{\^e}me
pour les $q$-diff{\'e}rences. Le deuxi{\`e}me r{\'e}sulte du fait
que l'{\'e}l{\'e}ment $g$ induit des changements de jauges
rationnels entre les syst{\`e}mes it{\'e}r{\'e}s par d{\'e}rivation
relatifs {\`a} $a$ et {\`a} $\overline{a}$.

\begin{remark}\label{remark:unic}

Supposons que $|q|<1$. Si on se restreint aux fonctions $a$ sans
z{\'e}ro ni p{\^o}le en $0$, il en sera de m{\^e}me de $\overline{a}$.
Si on impose alors au support du diviseur de $\overline{a}$
d'appartenir {\`a} la couronne ${\cal C} =\lbrace z \in
\mathbb{C}^*,  q <  |z| \leq 1 \rbrace$, la d{\'e}composition de $a$
sous forme standard est unique. En revanche si $0$ appartient au
support du diviseur de $a$, la relation
$\frac{\sigma_q(z^n)}{z^n}=q^n$ montre que $\overline{a}$ ne peut
{\^e}tre d{\'e}fini  qu'{\`a} multiplication pr{\`e}s  par un
{\'e}l{\'e}ment de $q^{\mathbb{Z}}$.
\end{remark}

\begin{lemme}\label{lemme:qnstand}
Soit $a \in K^*$. S'il existe un entier $n \neq 0$, un
{\'e}l{\'e}ment $h \in K^*$ et un nombre complexe $\mu \in
\mathbb{C}^*$ tels que $a^n= \mu  \frac{\sigma_q h}{h} $, alors il
existe  un {\'e}l{\'e}ment $g \in K^*$ et un nombre complexe
$\lambda \in \mathbb{C}^*$ tels que $a = \lambda \frac{\sigma_q g
}{g} $ et $\lambda^n \in \mu q^{\mathbb{Z}}$.
\end{lemme}

\textbf{D{\'e}monstration}\\
Il suffit d'{\'e}crire les factorisations de $h$ et de $a$  en
mon{\^o}mes, la relation sur $a^n$,  de comparer les deux formes
standards de $a^n$ auxquelles on aboutit et d'appliquer enfin la
remarque \ref{remark:unic}.

\begin{definition} Soient $a \in \mathbb{C}(z)^*$. On appelle
diviseur elliptique de $a$ l'image $div_E(a)$  de la partie premi{\`e}re
{\`a} $0$ du diviseur de $a$  par l'application naturelle de
$Div(\mathbb{C}^*)$ dans $Div(E)=Div(\mathbb{C}^*/q^{\mathbb{Z}})$.
\end{definition}
Autrement dit, si $a(z) = \lambda z^r \prod_{\alpha \in
\mathbb{C}^*} (z - \alpha)^{n_\alpha}$, on a $div_E(a) =
\sum_{\overline{\alpha} \in \mathbb{C}^*/q^{\mathbb{Z}}}
(\sum_{\alpha \equiv \overline{\alpha}}
n_\alpha).(\overline{\alpha})$. En {\'e}crivant $a$ sous forme
standard (lemme \ref{lemme:qstand}), on voit que :
\begin{lemme}\label{lemme:divq}
Soient $a \in \mathbb{C}(z)^*$. Alors $div_E(a)=0$ si et seulement
si il existe  un entier $r$, un nombre complexe $\mu \in \mathbb{C}^*$
 et $h \in K^*$ tel que $a=\mu z^r\frac{\sigma_q(h)}{h}$.
 \end{lemme}

\subsubsection{Solutions alg{\'e}briques}

\begin{prop}\label{prop:trans}
Soit $a \in K^*$ et
 soit $f \in \calM er(\mathbb{C}^*)$ une solution non nulle de l'{\'e}quation
 $\sigma_q (y) =a y (*)$. Alors,  $f$ est alg{\'e}brique sur $K_E$
si et seulement s'il existe  $\lambda  \in \mathbb{C}^*$ d'ordre
fini dans $\mathbb{C}^*/q^{\mathbb{Z}}$ et $g \in K^*$ tels que
$a(z) = \lambda \frac{\sigma_q(g)}{g}$.
\end{prop}
\textbf{D{\'e}monstration}\\
Consid{\'e}rons  le $\mathcal{D}_{K}$-module $\mathcal{A}$ de rang 1 associ{\'e} {\`a}
l'{\'e}quation $(*)$,  et le $\mathcal{D}_{K_E}$-module  $\mathcal{A}_E =  \mathcal{A} \otimes K_E$. Soit $G = Gal_{\omega_E}(\mathcal{A})$ le groupe de Galois de $\mathcal{A}_E$ relatif {\`a} $\omega_E$, et $\rho: G \rightarrow \bold{G}_m/_{C_E}$ la repr{\'e}sentation
correspondante, de degr{\'e} $1$, de $G$. Si
$\rho(G)$ est un sous groupe propre de $\bold{G}_m$,
alors il existe une entier $n \neq 0 $ tel que
$$\rho^{\otimes n}(G) \simeq  Gal_{\omega_E}(\mathcal{A}^{\otimes
  n})= \lbrace 1 \rbrace$$
Par {\'e}quivalence de cat{\'e}gories, le
$\mathcal{D}_{K_E}$-module $\mathcal{A}_E^{\otimes
  n}$ est donc trivial sur $K_E$ (voir aussi \cite{And} lemme $3.2.1.4$). D'apr{\`e}s le lemme \ref{lemme:comphom}, $\mathcal{A}^{\otimes
  n}$ est alors  trivial sur $K$, et il existe  $h \in K^*$ tel que
  \begin{equation}\label{eqn:tr1q}   a^{n}
  =\frac{\sigma_q(h)}{h}. \end{equation}
Le  lemme \ref{lemme:qnstand} permet alors de conclure.\\
R{\'e}ciproquement, s'il existe un tel couple $(\lambda \in
\mathbb{C}^*, g \in K^*)$, o\`u $\lambda$ est li\'e \`a $q$ par une
relation non triviale $\lambda^n = q^r$, alors $a^n
=q^r\frac{\sigma_q(g^n)}{g^n}=\frac{\sigma_q(g^n z^r)}{g^n z^r}$, et
$\mathcal{A}^{\otimes n}$  est trivial  dans $Diff(K, \sigma)$. Les
solutions de
  $\sigma_q(y)=ay$ v{\'e}rifient donc $y^n \in K.C_E$, et  sont bien  alg{\'e}briques sur $K_E$.

\subsubsection{Hypertranscendance des solutions d'{\'e}quations d'ordre $1$}

Le th{\'e}or{\`e}me suivant est un cas particulier du th{\'e}or{\`e}me
\ref{theorem:gen}, ou plus exactement, de la version plus pr{\'e}cise qu'il en sera donn{\'e}  plus bas (th{\'e}or{\`e}me \ref{theorem:gen2}), jointe {\`a} la remarque \ref{remark:gen2}.  Nous en d{\'e}taillons n{\'e}anmoins la d{\'e}monstration car
plusieurs de ses arguments  seront repris dans la preuve du th{\'e}or{\`e}me
g{\'e}n{\'e}ral.\\

\begin{theorem}\label{theorem:hyper}
Soit $a$  un {\'e}l{\'e}ment de $K^*=\mathbb{C}(z)^*$  et $f$ une
solution non nulle  dans $\mathcal{M}er(\mathbb{C}^*)$ de
l'{\'e}quation aux $q$-diff{\'e}rences $\sigma_q(f)=a f $. Alors
\begin{enumerate}

\item  $f$ et  $\partial f $ sont alg{\'e}briquement d{\'e}pendantes sur
  $C_E(z)$ si et seulement si $a$ est  de la forme $\mu
  \sigma_q(g)/g$, o{\`u} $\mu \in \mathbb{C}^*$ et $g \in \mathbb{C}(z)^*$.
\item $f$, $\partial f $ et ${\partial}^2f$ sont alg{\'e}briquement d{\'e}pendantes sur
  $C_E(z)$ si et seulement si $a$ est  de la forme $\mu z^r
  \sigma_q(g)/g$, $\mu  \in \mathbb{C}^*$, $r\ \in \mathbb{Z}$
  et $g \  \in \mathbb{C}(z)^*$.
\item Dans les autres cas, $f$ est hypertranscendante sur $C_E(z)$.
\end{enumerate}

\end{theorem}

Si $f$ est alg{\'e}brique sur $K_E$, la proposition \ref{prop:trans} montre que $a$ satisfait bien la condition du point 1) (avec $\mu$ d'ordre fini dans  $\mathbb{C}^*/q^{\mathbb{Z}}$). Sans perte de g{\'e}n{\'e}ralit{\'e},
on peut donc d{\'e}sormais supposer que $f$ est transcendante sur $K_E$. \\

Soit $n $ un entier $\geq 1$. Notons $t_n$  le degr{\'e} de transcendance du corps $L_n = K_E(f, {\partial}
 f ,...,\partial^{n}f)$ sur $K_E$, et $\delta_n$ la dimension du
 $\mathbb{C}$-sous-espace vectoriel  engendr{\'e} par  $
\frac{\partial a}{a},  \partial(\frac{\partial a}{a}), ..., \partial^{n-1}(\frac{\partial a}{a})$   dans  $K/(\sigma_q-1)(K)$.

\begin{lemme}\label{lemme:degtrq}
On suppose  $f$ transcendante sur $K_E$. Alors, $t_n = \delta_n + 1$.

\end{lemme}

\textbf{D{\'e}monstration}\\
On se r{\'e}f{\`e}re aux notations \ref{Not:ext} pour  la d{\'e}finition
de
$K_E(\mathcal{M})$.\\
Soit $f$  une solution non nulle de  l'{\'e}quation
$\sigma_q(f)=af$, $\mathcal{A}$ le $\mathcal{D}_K$-module
correspondant, et  $\mathcal{M}(n)$ l'objet de $Diff(K, \sigma_q)$
de repr{\'e}sentation matricielle $M(n)$
\begin{equation}\label{eqn:ext}
\left(\begin{array}{ccccccc}\\
        a & \cdots & \cdots & C_n^{k} \partial^{k}a & \cdots  &C_n^1 \partial^{n-1} a & \partial^{n} a   \\
        0 & a  & \cdots & \vdots    & \cdots  & \cdots   & \vdots\\
        0 &\cdots &  \ddots  &    C^{k-r}_{n-r}\partial^{k-r}a  &... &\cdots   & \partial^{n-r}a \\
        0   & ... &... & \cdots & \cdots & \cdots & \vdots\\
        0   & ... &... & a & \cdots & \cdots & \vdots\\
   0&...& ...& ... & \ddots&  \cdots & \partial^2 a \\
     0&...& ...& ... &...&  a & \partial a \\
         0& ...&...&...&...& 0 & a \end{array} \right).\end{equation}
         Une matrice fondamentale de solutions $U(n)$ en est donn{\'e}e par
\begin{equation}\label{eqn:extsol} \left(\begin{array}{ccccccc}\\
        f & \cdots & \cdots & C_n^{k} \partial^{k}f & \cdots  &C_n^1 \partial^{n-1} f & \partial^{n} f   \\
        0 & f  & \cdots & \vdots    & \cdots  & \cdots   & \vdots\\
        0 &\cdots &  \ddots  &    C^{k-r}_{n-r}\partial^{k-r}f  &... &\cdots   & \partial^{n-r}f \\
        0   & ... &... & \cdots & \cdots & \cdots & \vdots\\
        0   & ... &... & f & \cdots & \cdots & \vdots\\
   0&...& ...& ... & \ddots&  \cdots & \partial^2 f \\
     0&...& ...& ... &...&  f & \partial f \\
         0& ...&...&...&...& 0 & f \end{array} \right) \end{equation}

En effet, pour tout couple d'entiers $(k,r)$, $0 \leq r \leq n-k, 0
\leq k \leq n $, on a
$\sigma_q(\partial^{n-k-r}f)=\sum_{l=0}^{n-k-r} C^l_{n-k-r} \partial^l a
\partial^{n-k-r-l}f$.\\
 De plus,
$$M(n)U(n)_{(r,n-k)}= \sum_{l=0}^{n-k-r} C^l_{n-r}\partial^{l}a
C^{n-k-r-l}_{n-l-r}\partial^{n-k-r-l}f.
$$
Par cons{\'e}quent, on a, pour tout couple d'entiers $(k,r)$ :
$$\sigma_q( C^{n-k-r}_{n-r}\partial^{n-k-r}f)=\sum_{l=0}^{n-k-r}C^{n-k-r}_{n-r}C^l_{n-k-r} \partial^l a
\partial^{n-k-r-l}f$$
 $$=\sum_{l=0}^{n-k-r} C^l_{n-r}\partial^{l}a
C^{n-k-r-l}_{n-r-l}\partial^{n-k-r-l} f= M(n)U(n)_{(r,n-k)}  $$
car $C^{n-k-r}_{n-r}C^l_{n-k-r}=C^l_{n-r} C^{n-k-r-l}_{n-r-l}$, pour
tout entier $l \leq n-k-r$.

Pour tout $i=0,...,n-1$, soit par ailleurs
$\mathcal{E}_i$ le $\mathcal{D}_K$-module de repr{\'e}sentation
matricielle $\left( \begin{array}{cc}
   1 & \partial^{i} (\frac{\partial a}{a}) \\
0& 1 \end{array}. \right)$ C'est un  {\'e}l{\'e}ment de
 $Ext^1_{Diff(K, \sigma_q)} (\bold{1},\bold{1})$, qui correspond {\`a} l'{\'e}quation $\sigma_q(z)-z=
\partial^i(\partial a /a)$, dont une solution est donn{\'e}e par $\partial^{i}(\frac{\partial f}{f}) $. En effet,
on a $\frac{\sigma_q f}{f} = a$, d'o{\`u} $\sigma_q (\frac{\partial f}{f}) - \frac{\partial f}{f} =
\frac{\partial a}{a}$,  et en d{\'e}rivant $i$ fois cette {\'e}quation:   $\sigma_q(\partial^i( \frac{\partial f}{f})) -
\partial^i( \frac{\partial f}{f}) = \partial^i( \frac{\partial a}{a})$.  Ainsi, pour tout $i=0,...,n-1$, on a  $K_E(\mathcal{E}_i)=K_E(\partial^{i}(\frac{\partial f}{f}) )$ \\

Pour tout  entier $k < n$, la relation
$$\frac{\partial^{k+1}(f)}{f}= \frac{\partial^{k}(\frac{\partial f}{f}.f) }{f}
= \sum_{j=0}^{k}C^j_{k} \partial^{j}( \frac{\partial f}{f})
\frac{\partial^{k-j} f}{f}$$
  et une r{\'e}currence simple entra{\^\i}nent
que le corps $K_u:=K_E(\bigoplus_{i=1}^n \mathcal{E}_i)= K_E( \frac{\partial
f}{f},...,\partial^{n-1}(\frac{\partial
        f}{f}))$ est  {\'e}gal au corps $K'_u:=K_E(\frac{\partial
f}{f},  \frac{\partial^2 f}{f},...,\frac {\partial^{n}f}{f})$. Ce
dernier  est le corps de d{\'e}finition d'une matrice fondamentale
        de solutions de $\mathcal{M}(n)
\otimes \mathcal{A}^*$. Remarquons que  $\mathcal{M}(n)
\otimes \mathcal{A}^*$ est une extension it{\'e}r{\'e}e de l'objet unit{\'e}, puisque $\mathcal{A}
\otimes \mathcal{A}^* \simeq {\bold 1}$, ou plus concr{\`e}tement,  puisque sa matrice repr{\'e}sentative s'obtient en divisant par $a$ chacun des coefficients de celle de $\mathcal{M}(n)$, et n'a donc que des $1$ sur la diagonale.\\

Montrons que $Gal_{\omega_E}( \mathcal{M}(n)) \simeq Gal_{\omega_E}(\mathcal{M}(n)\otimes \mathcal{A}^*) \times Gal_{\omega_E}(
        \mathcal{A})$. Tout d'abord, les cat{\'e}gories tannakiennes $<\mathcal{M}(n)>$ et
$<\mathcal{M}(n) \otimes \mathcal{A}^* \oplus \mathcal{A}>$
co{\"\i}ncident. En effet, $\mathcal{A}$ {\'e}tant un sous objet de
$\mathcal{M}(n)$ l'une des  inclusions est triviale; l'autre
r{\'e}sulte du fait que $\mathcal{A}$,  de rang $1$,  est
trivialis{\'e} par tensorisation par son dual. Puisque la cat{\'e}gorie $< \mathcal{M}(n)\otimes \mathcal{A}^* \oplus
{\cal A}>$ est engendr{\'e}e par $\mathcal{M}(n) \otimes
\mathcal{A}^*$ et par ${\cal A}$, le groupe de Galois
$Gal_{\omega_E}(\mathcal{M}(n)\otimes \mathcal{A}^* \oplus
\mathcal{A})$
        est un sous groupe du produit direct
        $Gal_{\omega_E}(\mathcal{M}(n)\otimes \mathcal{A}^*) \times Gal_{\omega_E}(
        \mathcal{A})$ qui s'envoie surjectivement sur chacun des  facteurs. Mais $\mathcal{M}(n)\otimes \mathcal{A}^*$ est une extension it{\'e}r{\'e}e de l'objet $\bold{1}$, donc $Gal_{\omega_E}(\mathcal{M}(n)\otimes \mathcal{A}^*)$
        est un $C_E$-groupe unipotent, tandis que  $Gal_{\omega_E}(
        \mathcal{A})$ est un $C_E$-groupe semi-simple. Ils n'ont donc pas de quotients non triviaux isomorphes. Par cons{\'e}quent, $Gal_{\omega_E}(\mathcal{M}(n)\otimes \mathcal{A}^* \oplus \mathcal{A})$ remplit  tout le produit direct $Gal_{\omega_E}(\mathcal{M}(n)\otimes \mathcal{A}^*) \times Gal_{\omega_E}(
        \mathcal{A})$.\\

Le corps $L_n= K_E(\mathcal{M}(n))$
         est  le corps de d{\'e}finition d'une matrice fondamentale
        de solutions {\`a} coefficients dans $F$ de l'objet $\mathcal{M}(n)_E$ de $Diff(K_E, \sigma_q)$.
        De la   proposition \ref{prop:degtr}, on d{\'e}duit donc
$$t_n := deg.tr_{K_E} L_n = dim_{C_E}Gal_{\omega_E}( \mathcal{M}(n))=dim_{C_E}Gal_{\omega_E}( \mathcal{M}(n)\otimes \mathcal{A}^*) +  dim_{C_E}Gal_{\omega_E} ({\cal A}).$$

Une nouvelle application de  la proposition \ref{prop:degtr} montre que :
$$ dim_{C_E}Gal_{\omega_E}(\mathcal{M}(n)\otimes
\mathcal{A}^*) =  deg.tr_{K_E}K'_u = deg.tr_{K_E} K_u =    dim_{C_E}Gal_{\omega_E} (\bigoplus \mathcal{E}_i).$$

 Par
cons{\'e}quent,
$$t_n  =
dim_{C_E}Gal_{\omega_E}( \bigoplus \mathcal{E}_i) +
dim_{C_E}Gal_{\omega_E}(
        \mathcal{A}).$$
Comme $f$ est  transcendante sur $K_E$, la dimension  de $
Gal_{\omega_E}(
        \mathcal{A})$ est {\'e}gale {\`a} $1$.
 D'apr{\`e}s la proposition \ref{prop:indext1}, la
dimension de $ Gal_{\omega_E}(\bigoplus \mathcal{E}_i)$ est
{\'e}gale {\`a} la dimension $\delta_n$ du
 $\mathbb{C}$-sous-espace vectoriel de $K$  engendr{\'e} par les {\'e}l{\'e}ments  $
 \partial^i(\frac{\partial a}{a})$  ($i=0,...,n-1$)  modulo
 $(\sigma_q-1)(K)$. Ceci conclut la d{\'e}monstration du lemme \ref{lemme:degtrq}.\\

\begin{remark}
La m{\'e}thode d{\'e}crite dans l'introduction de l'article est ici
simplifi{\'e}e. En rang $1$, l'existence d'un isomorphisme entre
$Ext^1_{Diff(K, \sigma)} (\mathcal{A},\mathcal{A})$ et
$Ext^1_{Diff(K, \sigma)} (\bold{1},\bold{1})$ permet de remplacer
l'extension it{\'e}r{\'e}e  ${\mathcal{M}(n)}$ par la somme
directe des extensions simples $\bigoplus \mathcal{E}_i$ et de l'objet
$\mathcal{A}$. On aboutit ainsi {\`a} une situation proche des th{\'e}or{\`e}mes
de Kolchin et d'Ostrowski en th{\'e}orie des {\'e}quations diff{\'e}rentielles lin{\'e}aires.
Le m{\^e}me ph{\'e}nom{\`e}ne se produit pour les sommes directes
d'objets de rang 1 (voir le paragraphe 4.3.2  ci-dessous), mais pas
dans le cas g{\'e}n{\'e}ral. En particulier, nos r\'esultats ne
recouvrent pas ceux d'Ishizaki [12], qui correspondraient ici au cas
o\`u $\cal A$ est une extension non triviale de $\bold{1}$ par un
objet de rang $1$.
\end{remark}

\begin{remark}
En revanche, on peut  d{\'e}duire le th{\'e}or{\`e}me
\ref{theorem:hyper} du th{\'e}or{\`e}me de K.~Ishizaki \cite{Ishi}.
Celui-ci entra{\^\i}ne en effet que pour $a(z)$ de la forme $a_2(z)
= \prod (1 -a_iz)^{\alpha_i}$ i.e tel que $div(a_2) \subset
\mathbb{C}^*$, les solutions $y_2$ (non rationnelles) de
$\sigma_q(y)=a_2y$ sont hypertranscendantes. Par ailleurs, si $a$
est de la forme $a_1(z)=\mu z^r$, les solutions $y_1$ de
$\sigma_q(y)=a_1y$ v{\'e}rifient $\partial( \frac {\partial
y_1}{y_1})=0$. Donc, pour $a=a_1a_2$, $a_2\neq1$ les
solutions de $\sigma_q(y)=ay$ sont encore hypertranscendantes.\\

En revanche, les m{\'e}thodes d'Ishizaki ne suffisent pas {\`a}
d{\'e}montrer le r{\'e}sultat plus g{\'e}n{\'e}ral fourni par le
th{\'e}or{\`e}me \ref{theorem:gen}.

\end{remark}

\textbf{D{\'e}monstration du th{\'e}or{\`e}me
\ref{theorem:hyper}.}\\
D'apr{\`e}s le lemme  \ref{lemme:qstand}, on peut {\'e}crire
$a$ sous la
  forme  $ a = \overline{a} \frac{\sigma(f)}{f}$, $f \in K^*$
et $\overline{a}$ standard et les groupes de Galois aux
diff{\'e}rences associ{\'e}s aux extensions it{\'e}r{\'e}es par
d{\'e}rivation d{\'e}duites de  $a$ et $\overline{a}$ sont
{\'e}gaux. On peut donc supposer que $a$ est standard.

Soit donc $a \in K$  standard. Il s'agit maintenant de montrer que
le cas $(1): t_1 = 1$, c'est-{\`a} -dire en vertu du lemme
\ref{lemme:degtrq}, $\delta_1 = 0 $ (resp. $(2): t_2 = 2$,
c'est-{\`a}-dire $\delta_2 = 1$)   se produit   si et seulement si
$a$ est de la forme $\mu$ (resp.  $\mu z^r $ avec $ r$ un entier non
nul), et que sinon, c'est-{\`a}-dire si $a$ un un z{\'e}ro ou un
p{\^o}le non nul, on a: $\delta_n = n$, et donc $t_n = n+1$, pour
tout $n \geq 1$.
\\

\textbf{[preuve du point 1]}\\
Supposons que $\delta_1 = 0$, autrement dit que $\partial(a)/a$ soit
un {\'e}l{\'e}ment de  $(\sigma_q -1) \mathbb{C}(z)$. Alors, il
existe $k\in K$, de d{\'e}composition en {\'e}l{\'e}ments simples $
k(z) = \sum_{n=0}^D b_nz^n + \sum_{i=1}^t \sum_{l=1}^{\gamma_i}
\frac{\nu_i^l}{(z-c_i)^l}$, tel que
\begin{equation}\label{eqn:q1} \frac{\partial a}{a}=
  \sigma_q(k)-k\end{equation}
 Ecrivons $a$   sous la forme $a(z)=\mu
\prod_{i=1}^r(z-a_i)^{\alpha_i} $, o{\`u}  $a_i \neq q^{\mathbb{Z}}
a_{i'}$ si $i \neq i'$, et $\alpha_i  \in \mathbb{Z}$,   et $\mu \in
\mathbb{C}^*$. Puisque $\partial= z\frac{d}{dz}$, l'{\'e}quation
(\ref{eqn:q1}) s'{\'e}crit:
\begin{equation}\label{eqn:q12} \sum_i \alpha_i  +\sum_i \frac{\alpha_ia_i}{(z-a_i)}
 = \sum_{n=1}^D(q^nb_n-b_n)z^n
+\sum_i \sum_l
\frac{-\nu_i^l}{(z-c_i)^l}+\frac{\nu_i^l/q^l}{(z-c_i/q)^l}\end{equation}

 Il en r{\'e}sulte que  $ b_n=0 $ pour tout $ n \geq 1$ et $\sum \alpha_i
 =0$, et que $0$ n'est pas p{\^o}le de $k$. On va montrer maintenant  que $k$ ne peut avoir de p{\^o}le d'ordre
sup{\'e}rieur ou {\'e}gal {\`a} $1$.

En effet, supposons que $k$ ait un p{\^o}le   $c_{i_0}$ (non nul)
d'ordre $p>0$, et  consid{\'e}rons
 l'entier relatif  $n_0$ maximal  tel que $q^{-n_0}c_{i_0}$ soit un p{\^o}le
  d'ordre au moins   $p$ de $k$. Comme $ q^{-(n_0+1)} c_{i_0}$ est un p{\^o}le
  de $\sigma_q(k)$ et que $q^{-(n_0+1)} c_{i_0}$ n'est pas p{\^o}le d'ordre
  $p$ de $\sigma_q(k) -k$ ($a$ est standard), on en d{\'e}duit que $ q^{-(n_0+1)} c_{i_0}$ doit appara{\^\i}tre comme
  p{\^o}le d'ordre au moins $p$ de $k$. Ceci est absurde
par maximalit{\'e} de $n_0$.

On en d{\'e}duit que $k \in \mathbb{C}$ et que $a = \mu$ est  une fonction constante (au sens diff{\'e}rentiel).\\

R{\'e}ciproquement si $a$ est  constante, alors $\partial(a)/a=0 \in
(\sigma_q -1) \mathbb{C}(z)$, et $\delta_1 = 0$.

\medskip

\textbf{[preuve du point 2]}\\
Supposons que $\delta_2 = 1$, i.e. qu'il existe un nombre complexe $\lambda$ et un
{\'e}l{\'e}ment $k \in K$ tels que \begin{equation}\label{eqn:q2}
\partial(\partial a/a)
  +(\lambda) \partial a/a =\sigma(k)-k \end{equation}
Avec  les notations de la d{\'e}monstration
pr{\'e}c{\'e}dente, l'{\'e}quation (\ref{eqn:q2}) s'{\'e}crit
\begin{equation}\label{eqn:q22}
\lambda \sum_i \alpha_i + (\lambda-1)\sum_i \frac{\alpha_i
a_i}{z-a_i} -\sum_i\frac{\alpha_i a_i^2}{(z-a_i)^2} =
\sum_i\sum_l(\frac{\nu_i^l/q^l}{(z-c_i/q)^l}-\frac{\nu_i^l}{(z-c_i)^l})
+\sum (q^n b_n -b_n)z^n \end{equation} En suivant la
d{\'e}monstration du point $1$, on en d{\'e}duit que:
\begin{enumerate}
\item[a)] $k \in \mathbb{C}$.
\item[b)] le support du diviseur de  $a$ est r{\'e}duit {\`a} $(0)$.
\item[c)] $\lambda \sum_i \alpha_i=0$.
\end{enumerate}
De $b)$, il  r{\'e}sulte que $a$ est de la forme $a(z)=\mu z^r$
(avec $r =\sum_i \alpha_i$ non nul, sans quoi $\delta_2=0$, et
donc $\lambda=0$).\\

R{\'e}ciproquement,  si  $a(z)=\mu z^r$,
 alors $\partial(\partial a/a)=\partial(r) =0 \in (\sigma_q -1)
\mathbb{C}(z)$, et $\delta_2 \leq 1 $.\\

\textbf{[Preuve du point 3 ]}
\begin{lemme}\label{lemme:indel} Soit $a \in \mathbb{C}(z)^*$ telle
  que $a$ soit  standard et poss{\`e}de un p{\^o}le ou un z{\'e}ro non nul.
Alors la famille $\lbrace
\partial a / a ,..., \partial^j(\partial a /a),... \rbrace$ est lin{\'e}airement ind{\'e}pendante sur
$\mathbb{C}$ modulo $(\sigma_q -1)(\mathbb{C}(z))$.
\end{lemme}

\textbf{D{\'e}monstration}\\
Supposons qu'il existe $\lambda_{0},...., \lambda_N$ dans
$\mathbb{C}$, $\lambda_N \neq 0$ et $k \in \mathbb{C}(z)$, de la
forme $k(z) =\sum_{n=0}^D b_nz^n + \sum_{i=1}^t
\sum_{l=1}^{\gamma_i} \frac{\nu_i^l}{(z-c_i)^l}$, tels
que :\\

\begin{equation}\label{eqn:1} \sum_{j=0}^N \lambda_j \partial^j(\partial a/a)  =
 \sigma_q(k)-k = \sum_n(q^n b_n(qz) -b_n)z^n+
 \sum_i\sum_l(\frac{q^{-l}.\nu_i^l}{(z-c_i/q)^l}-\frac{\nu_i^l}{(z-c_i)^l}).\end{equation}
Une r{\'e}currence ais{\'e}e montre que pour tout entier $j \geq 0$
$$   \partial^j(\partial a /a) = \sum_i\frac{
  \alpha_i a_i^j(-1)^{j+1} j!}{(z-a_i)^{j+1}} + \ \mbox{des termes
  polaires d'ordre } \leq  j .$$

Soit alors $i_0$ tel que $a_{i_0} \neq 0$. Compte tenu des
  {\'e}critures pr{\'e}c{\'e}dentes et de l'unicit{\'e} de la d{\'e}composition en
  {\'e}l{\'e}ments simples sur $\mathbb{C}(z)$, on en d{\'e}duit que $a_{i_0}$ doit
  {\^e}tre un p{\^o}le d'ordre  $N+1$ de $\sigma_q (k) -k$, c'est-{\`a}-dire
  que soit $a_{i_0}$ soit $q a_{i_0}$ est un p{\^o}le d'ordre au
  au moins  $N+1$ de $k$.\\
Consid{\'e}rons  l'entier relatif  $n_0$ maximal  tel que $q^{-n_0}a_{i_0}$ soit un p{\^o}le
  d'ordre au moins  $N+1$ de $k$. Comme $ q^{-(n_0+1)} a_{i_0}$ est un p{\^o}le
  de $\sigma_q(k)$ et que $q^{-(n_0+1)} a_{i_0}$ n'est pas p{\^o}le d'ordre
  $N+1$ de $\sigma_q(k) -k$ ($a$ est standard), on en d{\'e}duit que $ q^{-(n_0+1)} a_{i_0}$ doit apparaitre comme
  p{\^o}le d'ordre au moins $N+1$ de $k$. Or
  ceci est absurde par maximalit{\'e} de $n_0$ et non nullit{\'e} de $a_{i_0}$.\\

Ainsi, pour tout entier $n \geq 1$, $\delta_n = n$  et $degtr_{K_E}
K_E(f,
\partial
 f ,..., \partial^n f...)=n+1$ ;  les fonctions  $f,\partial f ,..., \partial^n f $ sont donc alg{\'e}briquement ind{\'e}pendantes sur $K_E$ et $f$ est bien hypertranscendante sur $K_E$. Ceci conclut la preuve du  th{\'e}or{\`e}me \ref{theorem:hyper}.

\paragraph{Exemples}\label{par:qex}

On rappelle ici la d{\'e}finition des fonctions classiques de la
th{\'e}orie des $q$-diff{\'e}rences, qui illustrent chacun des cas
du th{\'e}or{\`e}me \ref{theorem:hyper}. On suppose que $|q| <1$.\\

On note  $\theta(z)$ (\cite{Ram2}) la fonction $\theta(z) =-\sum_{n\in \mathbb{Z}} (-1)^n
q^{\frac{n(n+1)}{2}}z^n$, qui correspond   {\`a} la fonction
\textit{theta} usuelle de la courbe elliptique $\mathbb{C}^* / q^{\mathbb{Z}}$. Elle est holomorphe sur $\mathbb{C}^*$, et
satisfait l'{\'e}quation aux $q$-diff{\'e}rences $\sigma_q(\theta)(z)= z^{-1}
\theta(z)$.\\

Pour $c\in \mathbb{C}^*$, on pose
$e_{q,c}(z)=\frac{\theta(z)}{\theta(z/c)}$ (voir (\cite{Sau}): c'est une solution m{\'e}romorphe sur  $\mathbb{C}^*$ de
l'{\'e}quation aux $q$-diff{\'e}rences $\sigma_q(e_{q,c})=c e_{q,c}$.\\

Un des analogues aux $q$-diff{\'e}rences de l'exponentielle est donn{\'e} par la fonction
$exp_q=\prod_{n\in \mathbb{N}^*}(1-(q-1)q^{-n}z)$, qui  est  m{\'e}romorphe sur
$\mathbb{C}$ et satisfait l'{\'e}quation aux $q$-diff{\'e}rences
$exp_q(qz)=(1-(q-1)z)exp_q (z)$.\\

\begin{enumerate}

\item [] \textit{Exemple 0} : $e_{q,-1}$ satisfait les hypoth{\`e}ses de
la proposition \ref{prop:trans} et v{\'e}rifie $e_{q,-1}^2 \in C_E$.
\item[] \textit{Exemple  1} : $e_{q,c}$ satisfait les hypoth{\`e}ses du point $1$
et v{\'e}rifie $\partial(e_{q,c})/e_{q,c} \in  C_E$
\item[]\textit{Exemple 2} :  Posons $l =\frac{\partial(\theta)}{\theta}$. La
  fonction $l$ correspond {\`a}  la fonction $\zeta$ de Weierstrass et
  v{\'e}rifie l'{\'e}quation aux $q$-diff{\'e}rences $\sigma_q(l)= l +1$. Par cons{\'e}quent, $\partial \frac{\partial \theta}{\theta} \in C_E$, et la
fonction $\theta$ v{\'e}rifie  une relation
alg{\'e}bro-diff{\'e}rentielle d'ordre 2 {\`a} coefficients dans $C_E(z)$

\item[]\textit{Exemple 3} : La fonction  $exp_q$ satisfait les hypoth{\`e}ses
  du point $3$. Elle est donc hypertranscendante sur $C_E(z)$.
\end{enumerate}

 \subsection{Syst{\`e}mes diagonaux aux $q$-diff{\'e}rences}

\subsubsection{Relations alg{\'e}briques entre solutions }

Soient $a_1 ,..., a_n$ $n$ {\'e}l{\'e}ments de $K^*$  et $A$ la
matrice diagonale de taille $n$ ayant comme coefficients diagonaux
les $a_i$. On consid{\`e}re le syst{\`e}me aux $q$-diff{\'e}rences
\begin{equation}\label{eqn:qder2}\sigma_q Y= \left( \begin{array}{cccc}
 a_1 & 0 & \cdots & 0 \\
 0 & a_2 & \cdots & 0\\
\vdots & 0& \ddots   & \vdots \\
0 & \cdots &0 & a_n  \\
\end{array} \right)    Y= A Y.\end{equation}

\begin{prop}\label{prop:qtransn} Soit $\Phi =\left( \begin{array}{cccc}
f_1 & 0 & \cdots & 0 \\
 0 & f_2 & \cdots & 0\\
\vdots & 0& \ddots   & \vdots \\
0 & \cdots &0 & f_n  \\
\end{array} \right)  \in Gl_n(\mathcal{M}er(\mathbb{C}^*))$ une matrice
 fondamentale de solutions de (\ref{eqn:qder2}).\\
 Alors, les fonctions  $f_i$ sont
  alg{\'e}briquement d{\'e}pendantes  sur $K_E$
si et seulement il existe des entiers $r_1, ...r_n$  non tous
 nuls,  et un {\'e}l{\'e}ment $h$ de $ K^*$  tels que
$\prod_i a_i^{r_i} =\frac{\sigma_q (h)}{h}$.
\end{prop}
(On notera que les diviseurs elliptiques des $a_i$ sont alors
lin{\'e}airement d{\'e}pendants sur $\mathbb{Z}$.)

\textbf{D{\'e}monstration}\\
Consid{\'e}rons  le $\mathcal{D}_{K}$-module  $\mathcal{A}$ de rang
$n$   associ{\'e} au syst{\`e}me (\ref{eqn:qder2}), et le
$\mathcal{D}_{K_E}$-module $\mathcal{A}_E=\mathcal{A}\otimes K_E$.
Notons, pour tout $i=1,...,n$, $\mathcal{A}_i$
 les  $\mathcal{D}_{K}$-module associ{\'e}s aux {\'e}quations
$\sigma_q(y)=a_i y$, de sorte  que $\mathcal{A} =\bigoplus_i
\mathcal{A}_i$.\\
 Soit $G=Gal_{\omega_E}(\mathcal{A})$ le groupe de Galois de $\mathcal{A}_E$ relatif {\`a} $\omega_E$, et
$\rho$ la repr{\'e}sentation correspondante de $G$ dans
 ${\bold{G}_m}^n/_{C_E}$.  Si $\rho(G)$ est un sous-groupe propre de
 ${\bold{G}_m}^n$, il est annul{\'e} par un caract{\`e}re $\chi$ non trivial, et   il existe des entiers $r_1 ,..., r_n$
 non tous nuls,  tels que
$$ \chi \circ \rho(G) \simeq  Gal_{\omega_E}(\bigotimes_i { \mathcal{A}_i }^{\otimes r_i}) =\lbrace 1
\rbrace.$$ Le
$\mathcal{D}_{K_E}$-module
$\bigotimes_i {(\mathcal{A}_i \otimes K_E)}^{\otimes r_i}$ est donc  trivial. D'apr{\`e}s le lemme \ref{lemme:comphom}, $\bigotimes_i
{\mathcal{A}_i}^{\otimes r_i}$ est alors  trivial sur $K$, et il
existe  $h \in K$ tel que :
 \begin{equation}\label{eqn:qtr} a_1^{r_1}...a_n^{r_n}
  =\frac{\sigma_q(h)}{h} \end{equation}

R{\'e}ciproquement, s' il existe  des entiers rationnels $r_1,...,r_n$  non
tous nuls
et un {\'e}l{\'e}ment  $h$ de $ K^*$ tels que  $ a_1^{r_1}...a_n^{r_n}
  =\frac{\sigma_q(h)}{h}$, alors $\bigotimes_i {\mathcal{A}_i}^{\otimes
    r_i}$   est   trivial dans $Diff(K, \sigma_q)$. Les
  solutions du syst{\`e}me (\ref{eqn:qder2}) v{\'e}rifient donc
  $y_1^{r_1}...y_n^{r_n} \in K.C_E$, et   sont bien  alg{\'e}briquement d{\'e}pendantes  sur $K_E$.

\subsubsection{Hyperind{\'e}pendance alg{\'e}brique}

Voici enfin, mis sous forme d'une condition n{\'e}cessaire et suffisante,
le th{\'e}or{\`e}me $1.1$ de l'introduction. On rappelle la d{\'e}finition $4.9$
des diviseurs elliptiques.

\begin{theorem}\label{theorem:gen2}
Soient $a_1 ,..., a_n$ des {\'e}l{\'e}ments non nuls  de $\mathbb{C}(z)$ et
$q$ un nombre complexe non nul de module diff{\'e}rent de $1$. Pour
tout $i=1 ,..., n$, soit $f_i \neq 0$  une  solution m{\'e}romorphe sur
$\mathbb{C}^*$  de l'{\'e}quation aux $q$-diff{\'e}rences
$\sigma_q(f_i)=a_i f_i$. Alors, les fonctions $f_1 ,..., f_n$ ainsi que leur d{\'e}riv{\'e}es
successives sont alg{\'e}briquement ind{\'e}pendantes sur $C_E(z)$   si et
seulement si les diviseurs elliptiques
des $a_i$ sont lin{\'e}airement ind{\'e}pendants sur $\mathbb{Z}$.
\end{theorem}
\begin{remark}\label{remark:gen2}
Si les diviseurs  elliptiques  sont li{\'e}s, on a d'apr{\`e}s le
lemme \ref{lemme:divq}, $\prod_i a_i^{r_i} = \mu z^r \sigma_q(
h)/(h) $ avec $r_i \in \mathbb{Z}, \mu \in \mathbb{C}, h \in K^*$.
En utilisant les fonctions classiques aux $q$-diff{\'e}rences
introduites au paragraphe \ref{par:qex}, on obtient $\prod_i
f_i^{r_i} = \lambda {\theta}^r  e_{q,\mu} h$, o{\`u} $\lambda \in
C_E$. En d{\'e}rivant une fois logarithmiquement, puis (si $r \neq
0$) en d{\'e}rivant une seconde fois, on obtient une relation de
d{\'e}pendance alg{\'e}brique  non triviale liant sur $K_E$ les
$f_i$ et leurs deriv{\'e}es. L'ordre
d'hyperalg{\'e}bricit{\'e} est encore une fois inf{\'e}rieur ou {\'e}gal {\`a} $2$, et  inf{\'e}rieur ou {\'e}gal {\`a} $1$ si $r = 0$.\\

\end{remark}

La preuve du th{\'e}or{\`e}me \ref{theorem:gen2} repose sur la g{\'e}n{\'e}ralisation  suivante du lemme \ref{lemme:degtrq}.

Soit $N $ un entier $\geq 1$. Notons $t_N$  le degr{\'e} de
transcendance du corps $L_N = K_E(\partial^{j}f_i ; i = 1,..., n; j = 0, ...N)$ sur $K_E$, et $\delta_N$ la dimension du
 $\mathbb{C}$-sous-espace vectoriel  engendr{\'e} par
les fonctions rationnelles $
\partial^{j}(\frac{\partial a_i}{a_i}), (i= 1,..., n; j = 0, ..., N-1)$   dans  $K/(\sigma_q-1)(K)$.

\begin{lemme}\label{lemme:degtrn}
Si les fonctions $f_1,...,f_n$ sont   alg{\'e}briquement ind{\'e}pendantes  sur $K_E$,
alors, $t_N = \delta_N + n$.
\end{lemme}

\textbf{D{\'e}monstration du lemme \ref{lemme:degtrn}}\\
Comme au paragraphe pr{\'e}c{\'e}dent, on note $A$ (resp. $\Phi$) la
matrice diagonale de coefficients $a_1, ..., a_n$ (resp. $f_1, ...,
f_n)$, et   $\mathcal{A}$ l'objet de  $Diff(K,\sigma_q)$  de
repr{\'e}sentation matricielle $A$, de sorte que $K_s := K_E(\mathcal{A}) = K_E(f_1, ..., f_n)$.\\

\medskip
Soit $\Phi$ une matrice fondamentatle de
        solutions de $\sigma_q(Y)=AY$ et  $\mathcal{M}(N)$ l'objet de $Diff(K, \sigma_q)$ de
repr{\'e}sentation matricielle

 \begin{equation}  \left(\begin{array}{ccccccc}\\
        A & \cdots & \cdots & C_m^{k} \partial^{k}A & \cdots  &C_m^1 \partial^{m-1} A & \partial^{m} A   \\
        0 & A  & \cdots & \vdots    & \cdots  & \cdots   & \vdots\\
        0 &\cdots &  \ddots  &    C^{k-r}_{m-r}\partial^{k-r}A  &... &\cdots   & \partial^{m-r}A \\
        0   & ... &... & \cdots & \cdots & \cdots & \vdots\\
        0   & ... &... & A & \cdots & \cdots & \vdots\\
   0&...& ...& ... & \ddots&  \cdots & \partial^2 A \\
     0&...& ...& ... &...&  A & \partial A \\
         0& ...&...&...&...& 0 & A \end{array} \right),  \end{equation}
dont une matrice fondamentale de solutions dans $F$ est donn{\'e}e
par :

\begin{equation}\label{eqn:extnsol} \left(\begin{array}{ccccccc}\\
        \Phi & \cdots & \cdots & C_m^{k} \partial^{k}\Phi & \cdots  &C_m^1 \partial^{m-1} \Phi & \partial^{m} \Phi   \\
        0 & \Phi  & \cdots & \vdots    & \cdots  & \cdots   & \vdots\\
        0 &\cdots &  \ddots  &    C^{k-r}_{m-r}\partial^{k-r}\Phi  &... &\cdots   & \partial^{m-r}\Phi \\
        0   & ... &... & \cdots & \cdots & \cdots & \vdots\\
        0   & ... &... & \Phi & \cdots & \cdots & \vdots\\
   0&...& ...& ... & \ddots&  \cdots & \partial^2 \Phi \\
     0&...& ...& ... &...&  \Phi & \partial \Phi \\
         0& ...&...&...&...& 0 & \Phi \end{array} \right) \end{equation}

 On a donc $L_N =  K_E(\partial^{j}f_i ; i = 1,..., n; j = 0, ...N)   = K_E(\mathcal{M}(N)  )$.

Pour  $i=1 ,..., n  \ j=0 ,..., N-1$, soit
   $\mathcal{E}_{i,j}$ l'extension de $\bold{1}$ par
$\bold{1}$ donn{\'e}e par le  syst{\`e}me  :

$$\mathcal{E}_{i,j}=\left (\begin{array}{cc}
1 & \partial^j (\frac{\partial a_i}{a_i}) \\
0& 1 \end{array} \right) $$ et soit
$\mathcal{E}(N)$ la somme directe $  \bigoplus \mathcal{E}_{i,j}$.
Ainsi, pour tout $i=1,...,n$ et tout $j=0,...,N-1$, on a
$K_E(\mathcal{E}_{i,j}) =K_E(\partial^j (\frac{\partial
f_i}{f_i}))$, et $K_E(\mathcal{E}(N))$ est le compositum de ces $Nn$
corps dans $F$.

Par r{\'e}currence sur l'entier $k <   N$, les relations
$$\frac{\partial^{k+1}(f_i)}{f_i}= \frac{\partial^{k}(\frac{\partial f_i}{f_i}) }{f_i}
= \sum_{j=0}^{k}C^j_{k} \partial^{j}( \frac{\partial f_i}{f_i})
\frac{\partial^{k-j} f_i}{f_i}$$ entra{\^\i}nent de nouveau que le
corps $K_u:=K_E(\mathcal{E}(N))= K_E(
\partial^{j}(\frac{\partial
        f_i}{f_i}); i = 1,..., n; j = 0, ..., N-1)$ co\" {\i}ncide avec le corps $K'_u:=K_E( \frac {\partial^{j}f_i}{f_i}; i = 1,..., n; j = 1, ..., N)$. Quant au corps $L_N$, c'est  clairement le compositum des corps $K_s$ et $K'_u$.

En d{\'e}finitive,  $L_N$ est le compositum des corps $K_s =
K_E(\mathcal{A})$ et $K_u = K_E(\mathcal{E}(N))$, d'o{\`u} $L_N =
K_E(\mathcal{A} \oplus \mathcal{E}(N)).$ De la proposition
\ref{prop:degtr}, on d{\'e}duit donc que
$$t_N  =  dim_{C_E}Gal_{\omega_E}( \mathcal{A} \oplus \mathcal{E}(N)).$$

Le groupe de Galois $Gal_{\omega_E}(\mathcal{E}(N) \oplus
\mathcal{A})$
        est un sous-groupe du produit direct
        $Gal_{\omega_E}(\mathcal{E}(N)) \times Gal_{\omega_E}(
        \mathcal{A})$ qui s'envoie surjectivement sur chacun des  facteurs. Comme  $\mathcal{E}(N) = \bigoplus_{i,j} \mathcal{E}_{i,j}$ est somme directe d'extensions de $\bold 1$ par $\bold 1$, son groupe de Galois $Gal_{\omega_E}(\mathcal{E}(N) )$
        est un $C_E$-groupe unipotent, tandis que $Gal_{\omega_E}(
        \mathcal{A})$ un $C_E$-groupe semi-simple. Par cons{\'e}quent,  $Gal_{\omega_E}(\mathcal{E}(N)  \oplus
\mathcal{A})$ remplit tout le produit direct
$Gal_{\omega_E}(\mathcal{E}(N))  \times Gal_{\omega_E}(
        \mathcal{A})$. Ainsi, \\
$$ t_N =dim_{C_E}Gal_{\omega_E} (\mathcal{A}) + dim_{C_E} Gal_{\omega_E}(\bigoplus \mathcal{E}_{i,j}) $$
Comme $f_1,...,f_n$ sont   alg{\'e}briquement ind{\'e}pendantes sur $K_E$, la dimension
de $ Gal_{\omega_E}(
        \mathcal{A})$ est {\'e}gale {\`a} $n$.
 D'apr{\`e}s la proposition \ref{prop:indext1}, la
dimension de $ Gal_{\omega_E}(\bigoplus \mathcal{E}_{i,j})$ est
{\'e}gale {\`a} la dimension $\delta_N$ du
 $\mathbb{C}$-sous-espace vectoriel de $K$  engendr{\'e} par les {\'e}l{\'e}ments  $
 \partial^j(\frac{\partial a_i}{a_i})$  ($i = 1, ..., n ; j=0, ..., N-1$)  modulo
 $(\sigma_q-1)(K)$. Ceci conclut la d{\'e}monstration du lemme \ref{lemme:degtrn}.\\

\bigskip

\textbf{D{\'e}monstration du th{\'e}or{\`e}me \ref{theorem:gen2}}\\
Pour {\'e}tablir le th{\'e}or{\`e}me  \ref{theorem:gen2}, il nous
reste, comme au point 3 de la preuve du th{\'e}or{\`e}me
\ref{theorem:hyper},  {\`a} v{\'e}rifier la g{\'e}n{\'e}ralisation
suivante du lemme \ref{lemme:indel}.

\begin{lemme}\label{lemme:derind}

Les fonctions
$\partial^j(\partial(a_i)/a_i)$  ($ i=1...n , j \in
\mathbb{N}$) sont lin{\'e}airement d{\'e}pendantes sur  $\mathbb{C}$  modulo $(\sigma_q-id)(\mathbb{C}(z))$ si et seulement si
 les diviseurs elliptiques $div_E(a_1),..., div_E(a_n)$  sont lin{\'e}airement d{\'e}pendants sur $\mathbb{Z}$.\\

\end{lemme}

\textbf{D{\'e}monstration}\\
Sans perte de g{\'e}n{\'e}ralit{\'e}, on peux supposer que les $a_i$
sont  standards. On fixe, une fois pour toutes,  une collection $S$
de repr{\'e}sentants de $\mathbb{C}^*/q^{\mathbb{Z}}$ dans
$\mathbb{C}^*$, not{\'e}e $S= \lbrace \alpha \in \mathbb{C}^*
\rbrace $. On peut alors {\'e}crire les $a_i$ sous la forme
 \begin{equation} \label{eqn:form} a_i= \mu_i z^{r_i} \prod_{\alpha \in
    S} (z - q^{n_{i,\alpha}}\alpha)^{\beta_{i,\alpha}} \end{equation}
avec  $\beta_{i,\alpha}, r_i \in
\mathbb{Z}$.\\

Supposons que les diviseurs elliptiques des  $a_i$ sont lin{\'e}airement
d{\'e}pendants sur $\mathbb{Z}$, et soient $l_1, ..., l_n$ des entiers non tous nuls tels que
\begin{equation}\label{eqn:liaisonq}
l_1 div_E(a_1) + l_2 div_E(a_2)+ ...+ l_n div_E(a_n)=0
\end{equation} dans $Div(\mathbb{C}^*/q^{\mathbb{Z}})$. On
d{\'e}duit de (\ref{eqn:liaisonq}) qu'il existe $r \in \mathbb{Z}$,
$\mu \in \mathbb{C}^*$ et $h \in K^*$ tels que : $\prod_i a_i^{l_i}
= \mu z^r \sigma_q( h)/(h)$. En d{\'e}rivant logarithmiquement
l'{\'e}quation pr{\'e}c{\'e}dente, et en d{\'e}rivant une seconde
fois on obtient :
$$ \sum_{i=1}^n l_i\partial(\frac {\partial a_i}{a_i}) =
\sigma_q(\partial(\frac{\partial h}{h}))-
\partial(\frac{\partial h}{h}).$$
On vient d'exhiber une liaison sur $\mathbb{C}$ entre les fonctions
$\partial^j(\partial(a_i)/a_i)$ $ i=1...n , j \in
\mathbb{N}$   (et m{\^e}me: $j = 0, 1, 2$) modulo $(\sigma_q-id)(\mathbb{C}(z))$.\\

\medskip

R{\'e}ciproquement, supposons qu'il existe une relation de d{\'e}pendance, qu'on peut choisir d'ordre $N$ minimal relativement  {\`a}  $j$,  liant les
$\partial^j(\partial(a_i)/a_i)$
  ($ i=1 ,...,  n, j \in
\mathbb{N}$) modulo $(\sigma_q-id)(\mathbb{C}(z)$ :
\begin{equation}\label{eqn:110q} \sum_{i=1}^n\sum_{j=0}^N \lambda^i_j
  \partial^j(\partial a_i/a_i)   =
 \sigma_q(f)-f , \end{equation}
 o{\`u} les $\lambda^i_j$ sont des nombres complexes, l'un des
 $\lambda^i_N$ est  non nul, et $f$ appartient {\`a} $\mathbb{C}(z)$.\\
On va montrer que, pour tout $\alpha \in S$ :
\begin{equation}\label{eqn:qliaison}
\sum_{i=1}^n \lambda^i_N \beta_{i, \alpha}=0
\end{equation}

On d{\'e}duit  des {\'e}quations  (\ref{eqn:qliaison}) que les
vecteurs  $\left( \begin{array}{c} \beta_{1, \alpha}\\
   \vdots \\
   \beta_{n, \alpha} \end{array} \right)$ , ${\alpha\in S}$, de $\mathbb{Z}^n$  sont li{\'e}s sur $\mathbb{C}$ et
donc sur $\mathbb{Z}$. Par cons{\'e}quent, ces {\'e}quations valent  avec des coefficients $\lambda^i_N$ entiers. Dans ces conditions,
 $$ \sum_{i=1}^n \lambda^i_N div_E(a_i)= \sum_{i=1}^n \lambda^i_N (\sum_{\alpha \in S}\beta_{i, \alpha}(\alpha))= \sum_{\alpha \in S}(\sum_{i=1}^n \lambda^i_N \beta_{i, \alpha})
 (\alpha)=0$$
\\
et les diviseurs elliptiques des $a_i$ sont bien lin{\'e}airement
d{\'e}pendants sur $\mathbb{Z}$.\\

Pour v{\'e}rifier (\ref{eqn:qliaison}),  {\'e}crivons la d{\'e}composition en {\'e}l{\'e}ments simples de $f$
sous  la forme $f(z) =\sum_{m=0}^M b_mz^m + \sum_{l=1}^p
\sum_{r=1}^{\gamma_l} \frac{\nu_l^r}{(z-d_l)^r}$, et reprenons les
expressions
$$ a_i= \mu_i z^{r_i}  \prod_{\alpha \in
    S} (z - q^{n_{i,\alpha}}\alpha)^{\beta_{i,\alpha}}$$
donn{\'e}es par la formule (\ref{eqn:form}); dans la suite, pour all{\'e}ger les notations, on posera
$b_{i,\alpha} = q^{n_{i,\alpha}}\alpha$ pour tout $i=1 ,...,  n, \alpha \in S$.\\

 Dans ces conditions, (\ref{eqn:110q}) s'{\'e}crit
\begin{equation}\label{eqn:11q} \sum_{i=1}^n\sum_{j=0}^N \lambda^i_j
  \partial^j(\partial a_i/a_i)   =
  \sum_m(q^m b_m(qz) -b_m)z^m+
 \sum_l\sum_r(\frac{q^{-r}.\nu_l^r}{(z-d_l/q)^r}-\frac{\nu_l^r}{(z-d_l)^r}).\end{equation}

Une r{\'e}currence ais{\'e}e montre que pour tout entier $ j \geq 0$
:
$$ \partial^j(\partial a_i /a_i) = \sum_{\alpha
  \in S} \frac{ \beta_{i,\alpha} (b_{i,\alpha})^j(-1)^{j+1} j!}{(z-b_{i,\alpha})^{j+1}} + \ \mbox{des termes
  polaires d'ordre } \leq  j.$$

Fixons un {\'e}l{\'e}ment  $\alpha$ de $S$, et  notons $I_{\alpha}$
l'ensemble d'indice
 $\lbrace \ i \in \{ 1, ..., n \} , \   \beta_{i,\alpha} \neq 0,\  \lambda^i_N
 \neq 0 \rbrace.$\\
 Tout d'abord, si $I_{\alpha}$ est vide, la relation
 (\ref{eqn:qliaison}) est v{\'e}rifi{\'e}e car, par d{\'e}finition de
 $I_{\alpha}$, on a alors $\beta_{i,\alpha} =0$ ou $\lambda^i_N=0$.\\

 En second lieu ,  la preuve du point $3$ de la d{\'e}monstration du
th{\'e}or{\`e}me \ref{theorem:hyper} assure que si $I_{\alpha}$ est
non vide, il ne peut {\^e}tre r{\'e}duit {\`a} un {\'e}l{\'e}ment.\\

Soit d{\'e}sormais $\alpha \in S$ tel que  $I_{\alpha}$ contient au
moins deux
  {\'e}l{\'e}ments. On note $ n_1< n_2 <...<n_t$ les valeurs
  distinctes, ordonn{\'e}es, prises par les $(n_{i,\alpha}), \ i \in I_{\alpha}$.\\
 Pour tout $l=1 ,...,  t$, on pose $I_{\alpha, n_l}=\lbrace \  i \in I_{\alpha}, \ \mbox{tels que
} \ n_{i,\alpha}= n_l \  \rbrace$.\\

La partie polaire d'ordre $N+1$ en la spirale $\alpha q^{\mathbb{Z}}$
du membre de gauche de l'{\'e}quation (\ref{eqn:11q}) est {\'e}gale
{\`a} $\sum_{l=1}^t \sum_{i \in I_{\alpha,n_l}}
\frac{\lambda^i_N\beta_{i,\alpha}(\alpha q^{n_l})^N(-1)^{N+1}
  N!}{(z - \alpha q^{n_l})^{N+1}}.$\\
On {\'e}crit la partie polaire d'ordre $N+1$ en la spirale $\alpha q^{\mathbb{Z}}$ de $f$ sous la forme $\sum_{k \in \mathbb{Z}}
\frac{\nu_k}{(z-\alpha q^k)}$, les $\nu_k$ {\'e}tant des nombres complexes presque
tous nuls.\\
En identifiant  la  partie polaire d'ordre $N+1$ en la spirale $\alpha q^{\mathbb{Z}}$ du
membre de gauche de l'{\'e}quation  {\`a} celle de $\sigma_q(f)-f$, on obtient
la relation suivante:
\begin{equation}\label{eqn:rec}
 \sum_{l=1}^t
\sum_{i \in I_{\alpha,n_l}} \frac{\beta_{i,\alpha}\lambda^i_N(\alpha
q^{n_l})^N(-1)^{N+1}
  N!}{(z - \alpha q^{n_l})^{N+1}} = \sum_{k \in \mathbb{Z}} \frac{\nu_{k+1}
  q^{-N-1} -\nu_k}{(z-q^k \alpha)^{N+1}}. \end{equation}
Par cons{\'e}quent,
\begin{enumerate}
\item[] \begin{equation}\label{eqn:2q} \forall m \notin \lbrace
    n_1 ,...,  n_t \rbrace, \ \nu_{m+1}q^{-N-1}-\nu_m=0. \end{equation}
\item[]  \begin{equation}\label{eqn:3q} \forall \   l = 1 ,...,  t,
  \  \nu_{n_l+1}q^{-N-1} -\nu_{n_l} = \sum_{i \in I_{\alpha,n_l}}\lambda^i_N\beta_{i,\alpha}(\alpha q^{n_l})^N(-1)^N
  N!.\end{equation}
\end{enumerate}
L'ensemble des entiers relatifs $m$ tels que $\nu_m$ soit non nul est
fini. La relation (\ref{eqn:2q}) entra{\^\i}ne que
\begin{enumerate}
\item $\nu_{n_1}=\nu_{n_{t +1}}=0.$
\item
  \begin{equation}
  \label{eqn:4q} \forall
    \ l=1 ,...,  t-1, \  \prod_{m = n_l +1}^{n_{l+1}-1} \frac{\nu_{m+1}}{\nu_m} =
    \prod_{m=n_l +1}^{n_{l+1}-1} q^{N+1} .\end{equation}

\end{enumerate}

De l'{\'e}quation (\ref{eqn:4q}) on tire les relations suivantes
$$ \forall l = 1,...,t-1,  \nu_{n_l+1} =\nu_{n_{l+1}} q^{(n_{l+1}-n_l -1)(N+1)}.$$ Alors, on a en reportant dans (\ref{eqn:3q}) que pour
tout $ l=1 ,...,  t-1$ :
\begin{equation}
\nu_{n_{l+1}}q^{-n_{l+1}(N+1)} -\nu_{n_l}q^{-n_l (N+1)} = \sum_{i \in
I_{\alpha,n_l}} \lambda^i_N\beta_{i,\alpha}{\alpha}^N(-1)^{N+1}
  N! \   \end{equation}
et \begin{equation}   -\nu_{n_t}q^{-n_t(N+1)} =\sum_{i \in
    I_{\alpha,n_t}}\lambda^i_N\beta_{i,\alpha}{\alpha}^N(-1)^{N+1} N!. \end{equation}\\

En sommant toutes ces {\'e}quations, on obtient :
$$\sum_{l=1}^t \sum_{i \in I_{\alpha,n_l}}\beta_{i,\alpha}\lambda^i_N \alpha^N (-1)^{N+1}
  N! =-\nu_{n_1}= 0,$$

d'o{\`u} $\sum_{i=1}^n \lambda^i_N \beta_{i, \alpha}=0$, puisque $I_{\alpha} = \bigcup_{l=1}^t I_{\alpha,n_l}$. Le lemme \ref{lemme:derind} est donc d{\'e}montr{\'e}.\\

\paragraph{Application}
Voici une illustration du th{\'e}or{\`e}me, pour l'analogue aux $q$-diff{\'e}rences de l'exponentielle (dont la d{\'e}finition a {\'e}t{\'e} rappel{\'e}e {\`a} fin du paragraphe $4.2$).

\begin{theorem}
Soient $\alpha_1 ,...,  \alpha_n$ des nombres complexes non nuls, qu'on suppose deux {\`a} deux distincts modulo $q^{\mathbb{Z}}$. Alors,
les fonctions $ exp_q(\alpha_iz), i=1 ,...,  n$ et leurs d{\'e}riv{\'e}es successives ne
v{\'e}rifient aucune relation alg{\'e}brique {\`a}
coefficients dans $K_E$.
\end{theorem}
Il suffit en effet d'appliquer le th{\'e}or{\`e}me \ref{theorem:gen2} avec $a_i = 1-(q-1)\alpha_i z$ pour tout $i = 1, .., n$, de sorte que $div_E(a_i) =(\frac{1}{(q-1)
 { \overline{\alpha_i}}})$, o{\`u}  $ \overline{ \alpha_i}$ d{\'e}signe la classe de $\alpha_i$ dans $C^*/q^Z$. Or les diviseurs $(\frac{1}{(q-1)
 {\overline{ \alpha_i}}})$ ne peuvent {\^e}tre lin{\'e}airement d{\'e}pendants sur $\mathbb{Z}$ que si deux d'entre eux au moins co{\"\i}ncident.

\section{Hypertranscendance des solutions d'{\'e}quations aux
  $\tau$-diff{\'e}rences }

\subsection{{\'E}nonc{\'e}s des r{\'e}sultats}

Nous montrons dans cette derni{\`e}re partie que les r{\'e}sultats pr{\'e}c{\'e}dents s'{\'e}tendent sans changement majeur (voir toutefois le `` point $A$" infra)  {\`a} l'{\'e}tude des $\tau$-diff{\'e}rences. Pr{\'e}cisons-en tout d'abord le cadre.

Soit $\tau  \in \mathbb{C}$ un nombre complexe non nul. On d{\'e}signe par $K=\mathbb{C}(z)$ le corps des fractions rationelles
{\`a} coefficients complexes, par $F=\mathcal{M}er(\mathbb{C})$ le corps
des fonctions m{\'e}romorphes sur $\mathbb{C}$ et par $\st$ l'automorphisme
de $F$ qui {\`a} $f(z) \in F$ associe $f(z+ \tau )$. On note $\Ct$
 le sous-corps de $F$ form{\'e} par les fonctions $\tau$-p{\'e}riodiques (c'est-{\`a}-dire les {\'e}l{\'e}ments fix{\'e}s par $\st$), et $\Kt=\Ct(z)$ le
compositum de $\Ct$
et de $K$ dans $F$. Les corps $K, F$ et $\Kt$ sont des corps aux diff{\'e}rences relativement {\`a} $\st$, admettant respectivement pour corps des $\st$-constantes $\mathbb{C},  \Ct$ et $\Ct$.\\

L'automorphisme $\st$ et la d{\'e}rivation $\partial = d/dz$
munissent le corps $F =  \calM er(\mathbb{C})$ et ses sous-corps
$\Kt$ et $K$ de structures de corps aux diff{\'e}rences
diff{\'e}rentiels, puisque $\st \partial = \partial \st$. Ces trois
corps admettent $\mathbb{C}$ comme corps de constantes
diff{\'e}rentielles.  Voici, dans ces conditions, les
$\tau$-analogues de la proposition \ref{prop:trans} et du th{\'e}or{\`e}me
\ref{theorem:hyper}.

\begin{prop}\label{prop:hypert}
Soient $a$ un {\'e}l{\'e}ment de $ K^*$ et $f   \in \mathcal{M}er(\mathbb{C})$ une solution non nulle  de
  l'{\'e}quation \begin{equation}\label{eqn:der2t} \st y=a y. \end{equation} Alors,

0. $f$ est alg{\'e}brique sur $\Kt$
si et seulement si $a$ est de la forme  $\zeta \frac{\st(g)}{g}$, o{\`u} $g \in K^*$ et $\zeta$  est une racine de l'unit{\'e} dans  $ \mathbb{C}^*$;

\begin{enumerate}

\item  $f$ et  $\partial f $ sont alg{\'e}briquement d{\'e}pendantes sur
  $\Ct(z)$ si et seulement si $a$ est  de la forme $\mu
  \sigma_q(g)/g$ o{\`u} $g \in \mathbb{C}(z)$ et  $\mu \in \mathbb{C}^*$;

\item Dans les autres cas, $f$ est hypertranscendante sur $\Ct(z)$.
\end{enumerate}

\end{prop}

\begin{definition} Soit $a \in \mathbb{C}(z)^*$.  On note $div_{\tau}(a)$
 l'image du diviseur  de $a$ par l'application naturelle de $Div(\mathbb{C})$ dans $Div(\mathbb{C}/\tau \mathbb{Z})$ et on l'appelle
diviseur p{\'e}riodique de $a$.\\
\end{definition}

Soit $A$ la matrice diagonale $(a_1, ..., a_n)$ : $$ \left(
\begin{array}{cccc}
a_1 & 0 &  \cdots   & 0 \\
 0 & a_2 &  \cdots   & 0\\
\vdots & 0& \ddots   & \vdots \\
0 &  \cdots   &0 & a_n  \\
\end{array} \right)  \in Gl_n(K).$$

Les $\tau$-analogues de la proposition \ref{prop:qtransn} et du
th{\'e}or{\`e}me \ref{theorem:gen2} s'{\'e}crivent :

\begin{theorem}\label{theorem:ttransn} Soit $\Phi =\left( \begin{array}{cccc}
f_1 & 0 &  \cdots   & 0 \\
 0 & f_2 &  \cdots   & 0\\
\vdots & 0& \ddots   & \vdots \\
0 &  \cdots   &0 & f_n  \\
\end{array} \right)  \in Gl_n(\mathcal{M}er(\mathbb{C}))$ une matrice
 fondamentale de solutions de $\st Y=A Y$. Alors,

 \begin{enumerate}

\item[i)] les fonctions  $f_1 ,...,  f_n$ sont
  alg{\'e}briquement d{\'e}pendantes  sur $\Kt$
si et seulement s'il existe des entiers $r_1,...,r_n$ non tous nuls
et un {\'e}l{\'e}ment $h \in K^*$ tels que
$a_1^{r_1}...a_n^{r_n}=\frac{\st(h)}{h}$.

\item [ii)] les fonctions  $f_1 ,...,  f_n$
v{\'e}rifient une relation alg{\'e}bro-diff{\'e}rentielle {\`a}
coefficients dans $\Kt$ si et seulement si les diviseurs
p{\'e}riodiques des $a_i$ sont $\mathbb{Z}$-lin{\'e}airement
d{\'e}pendants.

\end{enumerate}
\end{theorem}

\subsection{D{\'e}monstrations}

 Les d{\'e}monstrations de ces {\'e}nonc{\'e}s  rel{\`e}vent des m{\^e}mes
proc{\'e}d{\'e}s que pour les $q$-diff{\'e}rences. On n'en
explicitera donc pas les d{\'e}tails.

$A)$ $\tau$-analogue de  \ref{lemme:qstand}.$1$, \ref{remark:unic},
\ref{lemme:qnstand} et  \ref{lemme:divq}.\\

\begin{definition}[et lemme]\label{def:standt}Soit $a \in K^* = {\mathbb{C}(z)}^*$. On dit que $a$ est standard
 (chap.2, p.~29, lemme 2.1 et 2.2) si, pour tout $c \in \mathbb{C}$
dans le diviseur de $a$, et tout $n \in \mathbb{Z}-\lbrace 0
\rbrace$, $ c + n\tau$ n'apparait pas dans le diviseur de $a$.
Alors, il existe un  couple $(g, \overline{a}$), avec $g \in K^*$ et
$\overline{a}$ standard, tel que  $a = \overline{a}
\frac{\st(g)}{g}$. Une telle
d{\'e}composition est dite forme standard de $a$.\\
 Soient $f$ (resp.
$\overline{f}$) une solution m{\'e}romorphe sur $\mathbb{C}^*$ de
$\st(f)=af$ (resp. $\st(\overline{f})=\overline{a}\overline{f}$).
Pour tout $n \geq 0$, les corps $K(\partial^i(f);0 \leq  i \leq n)$
et $K(\partial^i(\overline{f}); 0 \leq  i \leq n)$ co{\"\i}ncident.
\end{definition}

\begin{lemme}\label{lemme:tstand}

 Soit $a$ un {\'e}l{\'e}ment de $K^*$.
\begin{enumerate}
\item Si on impose au support du diviseur de $\overline{a}$ d'appartenir {\`a} la
bande $C= \lbrace  z \in \mathbb{C}, \ 0 <
\mathcal{R}e(\frac{z}{\tau}) \leq 1 \rbrace$, la d{\'e}composition de
$a$ sous forme standard est alors unique.

\item S'il existe un entier $n$, un {\'e}l{\'e}ment $h
\in K^*$ et un nombre complexe $\mu \in \mathbb{C}^*$ tels que $a^n=
\mu \frac{\st h}{h}$, alors il existe  un {\'e}l{\'e}ment $g \in
K^*$ et une racine $n$-i{\`e}me de
l'unit{\'e} $\lambda \in \mathbb{C}^*$ tels que  $a =  \lambda
\frac{\st g }{g} $.

\item  $div_{\tau}(a)=0$ si et
seulement s'il existe un {\'e}l{\'e}ment $h \in K^*$ et  un nombre complexe $\mu \in \mathbb{C}^*$
  tels que $a=\mu \frac{\st(h)}{h}$.
\end{enumerate}
 \end{lemme}

La raison des diff{\'e}rences est que le point $(0)$ n'est plus,
contrairement aux $q$-diff{\'e}rences, fix{\'e} par l'op{\'e}rateur
$\sigma_{\tau}$. Ce ph{\'e}nom{\`e}ne se refl{\`e}te dans le fait
que la partie polaire en $0$ de $a$ n'intervient plus dans les
crit{\`e}res du paragraphe $5.1$. Dans le m{\^e}me esprit,  dans le
cadre des $\tau$-diff{\'e}rences l'ordre d'hyperalg{\'e}bricit{\'e}
des solutions est inf{\'e}rieur ou {\'e}gal {\`a} $1$, contrairement
au cas des $q$-diff{\'e}rences o{\`u} pour obtenir une relation
hyperalg{\'e}brique, on peut {\^e}tre amen{\'e} {\`a} consid{\'e}rer
les d{\'e}riv{\'e}es
secondes.\\

$B)$ V{\'e}rification des hypoth{\`e}ses \ref{hyp:ext}\\

Il s'agit de v{\'e}rifier ces hypoth{\`e}ses  pour le couple de
corps aux diff{\'e}rences form{\'e} de $K=\mathbb{C}(z)$ et de son
extension  $K' = \Kt = \Ct(z)$. On note   $G_{\tau} =
\Aut(\Ct/\mathbb{C}) $ (resp.  $G_K=Aut(\Kt/K)$) le groupe des
automorphismes  de l'extension $\Ct/\mathbb{C}$ (resp. $\Kt/K$), et
$\Ct(X)$ le corps des fractions
  rationnelles {\`a} coefficients dans  $\Ct$.  Pour tout  {\'e}l{\'e}ment  $u$ de $\mathbb{C} / \tau \mathbb{Z}$, la translation par $u$ d{\'e}finit un automorphisme du corps $\Ct$, dont on note $\gamma_u$ le prolongement canonique {\`a} $\Ct(X)$, d{\'e}fini par son action sur les coefficients.  L'ensemble $\Gamma = \lbrace \gamma_u, u \in \mathbb{C} / \tau
\mathbb{Z} \rbrace$ forme un sous-groupe de $G_K$, isomorphe {\`a}
$\mathbb{C} / \tau \mathbb{Z}$. Le lemme $4.2$ est alors
inchang{\'e}, le seul point {\`a} v{\'e}rifier est que  l'analogue
aux $\tau$-diff{\'e}rences du lemme de \cite{Eti} est toujours
valable. Ce fait, classique,  repose sur la transcendance de la
fonction $z$ sur le corps $\Ct$ et {\'e}nonce  :
\begin{lemme}\label{lemme:teti}
Soient $(c_0,...,c_N)$ une famille d'{\'e}l{\'e}ments de $\Ct$. On
suppose  que :
\begin{equation}\label{eqn:teti} \sum_{i=0}^N
c_i(z)z^i=0.\end{equation} Alors, pour tout $i=0,...,N$, $c_i=0$.
\end{lemme}

\noindent

$C)$ $\tau$-analogues des lemmes $4.16$ et
\ref{lemme:derind}.\\
Il suffit de supprimer \textit{non nul} dans le lemme 4.16, et, plus
g\'en\'eralement, de remplacer $div_E$ par $div_{\tau}$ dans le
lemme 4.16, puisque le point $0$ ne joue plus
de r{\^o}le particulier (voir le lemme \ref{lemme:tstand}  et le  commentaire qui le suit).\\

\textbf{Application}\\
Dans cette application $\tau=1$, et on fixe des nombres complexes
$\alpha_1,...,\alpha_k$ ainsi que deux entiers $n,m \geq 2$. La
relation de distribution $$\Gamma(z) \Gamma(z
+\frac{1}{n})...\Gamma(z+\frac{n-1}{n})=(2
\pi)^{\frac{n-1}{2}}n^{1/2-nz} \Gamma(nz)$$ (et l'{\'e}quation
diff{\'e}rentielle $\frac{d}{dz}m^z=(logm)m^z$) montrent que si $k=n$ et
si $\alpha_i$ est congru {\`a} $\frac{i}{n}$ modulo $\mathbb{Z}$
pour tout $i$, alors les fonctions
$\Gamma(z+\alpha_1),...,\Gamma(z+\alpha_n),\Gamma(nz)$ sont
hyperalg{\'e}briquement d{\'e}pendantes sur $K$ ; si de plus $m$ et
$n$ sont multiplicativement d{\'e}pendants, alors ces fonctions
ainsi que $m^z$ sont alg{\'e}briquement d{\'e}pendantes sur $K$. Inversement :

 \begin{coro}

 Soient $\alpha_1  ,...,   \alpha_k$ des {\'e}l{\'e}ments de
$\mathbb{C}$ et  $ n_1, ..., n_h, m$ des entiers $\geq 2$. On
suppose que les diviseurs p{\'e}riodiques $( \overline{\alpha_i}),i=
1, ..., k, \ (\sum_{i=0}^{n_j-1}(\overline{\frac{i}{n_j}}), j = 1,
..., h$ sont $\mathbb{Z}$-lin{\'e}airement ind{\'e}pendants. Alors :
\begin{enumerate}
\item   les fonctions $\Gamma(z+\alpha_1),...,\Gamma(z+\alpha_k),\Gamma(n_1z),..., \Gamma(n_hz)$ et $m^z$ son alg{\'e}briquement
ind{\'e}pendantes
 sur $K_{\tau}$;
\item les
fonctions $\Gamma(z +\alpha_i), \Gamma(n_jz), i= 1, ..., k, j = 1,
...h$ sont hyperalg{\'e}briquement ind{\'e}pendantes sur $\Kt$.
\end{enumerate}
\end{coro}

\textbf{D{\'e}monstration}\\
On a :
\begin{enumerate}
\item[] $\st(\Gamma(z +\alpha_i)= a_i(z)\Gamma(z
  +\alpha_i)$, avec $a_i(z)=z+\alpha_i$
\item[] $\st(\Gamma(n_jz))= a'_j(z) \Gamma(nz)$, avec
$a'_j (z)=\prod_{i=1}^{n_j}(n_j z + i)$
\item[] $\st(m^z)=a_0(z) m^z$, avec $a_0(z)=m$.
\end{enumerate}

\medskip

$1)$ D'apr{\`e}s le  th{\'e}or{\`e}me \ref{theorem:ttransn}, point
$i)$, les fonctions $\Gamma(z +\alpha_i)$, $\Gamma(n_jz)$, $m^z$
sont  alg{\'e}briquement d{\'e}pendantes  sur $\Kt$  si et seulement s'il
existe des entiers $r_1,...,r_k,l_{1},..., l_{h},r $ non tous
nuls et un {\'e}l{\'e}ment $h \in K^*$ tels que $$m^r \prod_{j=1}^h
\prod_{i=0}^{n_j-1}(n_jz +i)^{l_j} \prod_{i=1}^k (z +\alpha_i)^{r_i}
=\frac{\st(h)}{h}.$$ Cette relation entra{\^\i}ne :
\begin{enumerate}
\item $\sum_{j=1}^h l_j (\sum_{i=0}^{n_j-1}(\overline{\frac{i}{n_j}}))
+ \sum_{i=1}^k r_i ( \overline{\alpha_i}) =0$;
\item $m^r\prod_{j=1}^h (n_j)^{l_j n_j}=1.$
\end{enumerate}

Il r{\'e}sulte de l'hypoth{\`e}se d'ind{\'e}pendance des diviseurs p{\'e}riodiques
que : $l_1=...=l_h=r_1=...=r_k=0$. Donc, $r\neq 0$ et $m^r=1$, ce
qui est absurde car $m \geq 2$. Les fonctions $\Gamma(z +\alpha_i)$,
$\Gamma(n_jz)$, $m^z$ sont donc  alg{\'e}briquement
ind{\'e}pendantes sur $\Kt$.\\

$2)$ D'apr{\`e}s le  th{\'e}or{\`e}me \ref{theorem:ttransn}, les
fonctions  $\Gamma(z +\alpha_i)$, $\Gamma(n_jz)$ ($i=1,...,k$,
$j=1,...,h$) sont hyperalg{\'e}briquement ind{\'e}pendantes sur
$\Kt$ si et seulement si
 les diviseurs p{\'e}riodiques $(- \overline{\alpha_i}),
\sum_{i=0}^{n_j -1}(-\overline{\frac{i}{n_j}}) $  sont
$\mathbb{Z}$-lin{\'e}airement ind{\'e}pendants. La conclusion du
point $2$ du corollaire en r{\'e}sulte. On peut aussi exprimer cette
conclusion de la fa{\c c}on suivante. Pour tout entier $t \geq 0$,
soit $\psi^{(t)}=(\frac{d}{dz})^{t+1} Log \Gamma$ la fonction
polygamma d'ordre $t$. Alors, les fonctions $\Gamma(z+\alpha_i)$,
$\psi^{t}(z+ \alpha_i)$, $\Gamma(n_j z)$, $\psi^{(t)}(n_j z)$
($i=1,...,k$, $j=1,...,h$ et $t \in \mathbb{N}$) sont
alg{\'e}briquement ind{\'e}pendantes sur $K_{\tau}$.

\end{document}